\documentclass[12pt]{amsart}
\usepackage{amssymb}

\begin{document}
\numberwithin{equation}{section}

\def\Label#1{\label{#1}}

\def\1#1{\ov{#1}}
\def\2#1{\widetilde{#1}}
\def\3#1{\mathcal{#1}}
\def\4#1{\widehat{#1}}

\def\s{s}
\def\k{\kappa}
\def\ov{\overline}
\def\span{\text{\rm span}}
\def\tr{\text{\rm tr}}
\def\GL{{\sf GL}}
\def\xo {{x_0}}
\def\Rk{\text{\rm Rk\,}}
\def\sg{\sigma}
\def \emxy{E_{(M,M')}(X,Y)}
\def \semxy{\scrE_{(M,M')}(X,Y)}
\def \jkxy {J^k(X,Y)}
\def \gkxy {G^k(X,Y)}
\def \exy {E(X,Y)}
\def \sexy{\scrE(X,Y)}
\def \hn {holomorphically nondegenerate}
\def\hyp{hypersurface}
\def\prt#1{{\partial \over\partial #1}}
\def\det{{\text{\rm det\,}}}
\def\wob{{w\over B(z)}}
\def\co{\chi_1}
\def\po{p_0}
\def\fb {\bar f}
\def\gb {\bar g}
\def\Fb {\ov F}
\def\Gb {\ov G}
\def\Hb {\ov H}
\def\zb {\bar z}
\def\wb {\bar w}
\def \qb {\bar Q}
\def \t {\tau}
\def\z{\chi}
\def\w{\tau}
\def\Z{\zeta}

\def \T {\theta}
\def \Th {\Theta}
\def \L {\Lambda}
\def\b{\beta}
\def\a{\alpha}
\def\o{\omega}
\def\l{\lambda}

\def \im{\text{\rm Im }}
\def \re{\text{\rm Re }}
\def \Char{\text{\rm Char }}
\def \supp{\text{\rm supp }}
\def \codim{\text{\rm codim }}
\def \Ht{\text{\rm ht }}
\def \Dt{\text{\rm dt }}
\def \hO{\widehat{\mathcal O}}
\def \cl{\text{\rm cl }}
\def \bR{\mathbb R}
\def \bC{\mathbb C}
\def \bP{\mathbb P}
\def \C{\mathbb C}
\def \bL{\mathbb L}
\def \bZ{\mathbb Z}
\def \bN{\mathbb N}
\def \scrF{\mathcal F}
\def \scrK{\mathcal K}
\def \scrM{\mathcal M}
\def \cR{\mathcal R}
\def \scrJ{\mathcal J}
\def \scrA{\mathcal A}
\def \scrO{\mathcal O}
\def \scrV{\mathcal V}
\def \scrL{\mathcal L}
\def \scrE{\mathcal E}
\def \hol{\text{\rm hol}}
\def \aut{\text{\rm aut}}
\def \Aut{\text{\rm Aut}}
\def \J{\text{\rm Jac}}
\def\jet#1#2{J^{#1}_{#2}}
\def\gp#1{G^{#1}}
\def\gpo{\gp {2k_0}_0}
\def\emmp {\scrF(M,p;M',p')}
\def\rk{\text{\rm rk}}
\def\Orb{\text{\rm Orb\,}}
\def\Exp{\text{\rm Exp\,}}
\def\ess{\text{\rm Ess\,}}
\def\mult{\text{\rm mult\,}}
\def\Jac{\text{\rm Jac\,}}
\def\Span{\text{\rm span\,}}
\def\d{\partial}
\def\D{\3J}
\def\pr{{\rm pr}}
\def\dbl{[\![}
\def\dbr{]\!]}
\def\nl{|\!|}
\def\nr{|\!|}

\def \D{\text{\rm Der}\,}
\def \Rk{\text{\rm Rk}\,}
\def \ima{\text{\rm im}\,}
\def \vfi{\varphi}

\title[Transversality of CR mappings]
{Transversality of CR mappings}
\author[P. Ebenfelt and L. P. Rothschild]{Peter Ebenfelt and
Linda P. Rothschild} \footnotetext{{\rm The first author is
supported in part by DMS-0100110. The second author is
supported in part by DMS-0100330.\newline}}
\address{ Department of Mathematics, University of California
at San Diego, La Jolla, CA 92093-0112, USA}
\email{pebenfel@math.ucsd.edu, lrothschild@ucsd.edu }

\date{\number\year-\number\month-\number\day}

\abstract We prove here new results about
transversality and related geometric properties of a 
holomorphic, formal, or CR
mapping, sending one generic submanifold of $\bC^N$
into another. One of our main results is that
a finite mapping is transversal to the
target manifold provided this manifold is of finite type.
For the case of hypersurfaces, transversality in this context 
was proved by Baouendi and the second author in 1990. The
general case of generic manifolds of higher codimension,
which we treat in this paper, had remained an open problem
since then. Applications of this result include a
sufficient condition for a finite mapping to be a local
diffeomorphism.
\endabstract

\newtheorem{Thm}{Theorem}[section]
\newtheorem{Def}[Thm]{Definition}
\newtheorem{Cor}[Thm]{Corollary}
\newtheorem{Pro}[Thm]{Proposition}
\newtheorem{Lem}[Thm]{Lemma}
\newtheorem{Rem}[Thm]{Remark}
\newtheorem{Ex}[Thm]{Example}

\maketitle
\section{Introduction}
The focus of this paper is to prove new results about
transversality and related geometric properties of a holomorphic
(or formal) mapping that sends one generic submanifold of
$\bC^N$ into another. Our main result is a sufficient
condition for such a mapping to be transversal to the target
manifold. We also give applications and refinements of this
result such as e.g.\ a criterion for a finite mapping to be a
local (formal) biholomorphism.

In order to state our results more precisely, we need to
introduce some notation. Recall that a real smooth submanifold $M$
of codimension $d$ through $p_0$ in $\bC^N$ is called {\it
generic} (at $p_0$) if there is a smooth function $\rho\colon
(\bC^N,p_0)\to\bR^d$ such that $M$ is locally given near $p_0$ by
$\rho=0$ and
$$
\partial\rho_1\wedge\ldots \wedge \partial\rho_d\neq 0
$$ at $p_0$. A generic submanifold $M$ is a CR manifold with
the CR structure defined by the CR bundle $T^{0,1}M:=\bC TM\cap
T^{0,1}\bC^N$ and its complex conjugate
$T^{1,0}M:=\ov{T^{0,1}M}$. Here, $\bC TM$ denotes the
complexified tangent space of $M$ and $T^{0,1}\bC^N$ denotes
the bundle of $(0,1)$ vectors on $\bC^N$. Let $\widetilde
M\subset \bC^N$ be a generic submanifold and $\tilde p_0 \in
\widetilde M$. A germ of a holomorphic mapping
$H\colon (\bC^N,p_0)\to (\bC^N,\tilde p_0)$ is said to be {\it
CR transversal} to 
$\widetilde M$  at $p_0$
if
\begin{equation}\Label{e:CRtransdef}
T^{1,0}_{\tilde p_0}\widetilde M+dH(T^{1,0}_{p_0}\bC^N)=T^{1,0}_{\tilde p_0}\bC^N.
\end{equation}
One significance of CR transversality is stability in the
following sense. If $H$ is CR transversal to $\widetilde
M$, then the inverse image
$H^{-1}(\widetilde M)$ is a germ of a generic submanifold and for
every sufficiently small smooth perturbation
$\widetilde M_{\epsilon}$ of $\widetilde M$, the inverse image
$H^{-1}(\widetilde M_\epsilon)$ is a smooth generic perturbation of
$H^{-1}(\widetilde M)$. Recall that if $f\colon(\bR^k,p_0) \to
(\bR^\ell,\tilde p_0)$ is a germ of a smooth map, and if
$\widetilde M
\subset \bR^\ell$ is a smooth manifold through $\tilde p_0$,
then
$f$ is called {\em transversal} to $\widetilde M$ at $p_0$ if 
\begin{equation}\Label{e:transdef}
T_{\tilde p_0}\widetilde M+df(T_{p_0}\bR^k)=T_{\tilde p_0}\bR^{\ell}
\end{equation}
 (see e.g.\ \cite{GG86}). It is not hard to show that if
$H$ is CR transversal to $\widetilde M$ at $p_0$, then $H$,
regarded as a smooth mapping
$(\bR^{2N},p_0)\to (\bR^{2N},\tilde p_0)$ is transversal to
$\widetilde M$ at $p_0$.
The converse is not true in general. However,
the  notions of transversality and CR transversality
coincide if $H^{-1}(\widetilde M)$ is generic at
$p_0$ (see Section \ref{s:trans}).

One of the main results in this paper, Theorem
\ref{t:main}, is a sufficient condition for a holomorphic
mapping
$H\colon (\bC^N,p_0)\to (\bC^N,\tilde p_0)$ to be CR
transversal to a generic submanifold $\widetilde M\subset \bC^N$
through $\tilde p_0$. Our condition is also necessary in
some cases  (see e.g. Theorem \ref{t:maincor2}).  A special case of
Theorem \ref{t:main} is the following. Recall that a generic
submanifold $M$ is said to be of {\it finite type}  at
$p_0$ (in the sense of Kohn, and Bloom--Graham) if the
(complex) Lie algebra
$\frak g_M$ generated by all smooth $(1,0)$ and $(0,1)$
vector fields tangent to $M$ satisfies $\frak
g_M(p_0)=\bC T_{p_0}M$.

\begin{Thm}\Label{t:maincor} Let $\widetilde M\subset\bC^N$
be a smooth generic submanifold of finite type at $\tilde
p_0$ and $H\colon (\bC^N,p_0)\to (\bC^N,\tilde p_0)$ a finite
holomorphic mapping. If there exists a smooth generic
submanifold $M\subset \bC^N$ through $p_0$ of the
same dimension as $\widetilde M$ with
$H(M)\subset \widetilde M$, then $H$ is CR transversal to
$\widetilde M$ at
$p_0$.
\end{Thm}

%

We would like to point out that since $H$ is assumed to be a
finite mapping, the condition that $\widetilde M$ is of finite type at
$\tilde p_0$ could be replaced by the condition that $M$ is of
finite type at $p_0$ (see Proposition \ref{p:ft}). In the case
where
$M$ and $\widetilde M$ are assumed to be hypersurfaces, 
Theorem
\ref{t:maincor} follows from results by
Baouendi and the second author (see \cite{BRgeom}). The case
of generic submanifolds of higher codimension had remained an
open problem since the publication of \cite{BRgeom}. 

Theorems
\ref{t:maincor} and
\ref{t:main} can be viewed as formal Hopf Lemma type results.
For smooth mappings between hypersurfaces where the target has
some convexity properties, results of this type have been known
for some time; see Fornaess \cite{Forn76},
\cite{Forn78} for the case of pseudoconvex hypersurfaces and 
\cite{BRhopf} for a more general convexity condition. See also
\cite{BHR95} by Baouendi, X. Huang, and the second author,
and \cite{HP96} by Huang and Y. Pan for other types of Hopf
lemmas for smooth mappings between hypersurfaces. Other
results along these lines in the higher codimensional case can 
be found in e.g.\  the papers \cite{CR94}, \cite{CR98} by
Chirka and Rea. 

The other results in this paper are refinements and
applications of Theorem \ref{t:main}. It is easy to find
examples (see Section \ref{s:new}) of finite holomorphic
mappings satisfying the hypotheses of Theorem
\ref{t:maincor} (and, hence, CR transversal) that are not
biholomorphic. In our next result, we give a sufficient
condition on the source manifold $M$ that guarantees that
any finite holomorphic mappings sending
$M$ into a generic submanifold is biholomorphic.  Recall
that $M$ is said to be {\it finitely nondegenerate} at $p_0$ if
\begin{equation}\Label{e:fnd}
\span_\bC\left \{ L^\alpha\left(\frac{\partial
\rho^j}{\partial Z}\right)(p_0)\colon j=1,\ldots, d,\
\alpha\in
\mathbb N_+^n\right
\}=\bC^N,
\end{equation} where $\span_\bC$ denotes the vector
space spanned over $\bC$ and
$L^\alpha:=L_1^{\alpha_1}\ldots L_n^{\alpha_n}$.  Here,
$L_1,\ldots, L_n$ is a basis for the smooth
$(0,1)$ (or CR) vector fields tangent to $M$ near $p_0$, and
$\rho=(\rho_1,\ldots,\rho_d)^t$ is a defining function for $M$
near $p_0$. The following is a consequence of our Theorem
\ref{t:fndsimple}.

\begin{Thm}\Label{t:maincor1} Let $M\subset\bC^N$ be a
smooth generic submanifold  of finite type and finitely
nondegenerate at $p_0\in M$. If $\widetilde M\subset\bC^N$ is a
smooth generic submanifold of the same dimension and $H\colon
(\bC^N,p_0)\to \bC^N$ is a germ of a finite holomorphic mapping
with
$H(M)\subset \widetilde M$, then $H$ is
a local biholomorphism at $p_0$.
\end{Thm}

An early result of this type is contained in a paper by
Pinchuk \cite{P77} in which he proves that nonconstant CR
mappings between strictly pseudoconvex hypersurfaces are local
diffeomorphisms. For self maps of hypersurfaces, see also
Bedford and Bell \cite{BB82} as well as 
\cite{BRhopf}.

For essentially finite generic manifolds of
finite type, the converse of Theorem \ref{t:maincor} is also
true. (See Section \ref{s:new}.) More
precisely, we have the following result.

\begin{Thm}\Label{t:maincor2} Let $M\subset\bC^N$ be a
smooth generic submanifold  of finite type and essentially
finite at
$p_0$. Suppose
$\widetilde M\subset\bC^N$ is  a smooth generic submanifold
of the same dimension, and  $H\colon (\bC^N,p_0)\to \bC^N$  
a germ of a holomorphic mapping with
$H(M)\subset \widetilde M$.  Then the following are equivalent.
\medskip

{\rm (i)} $H$ is CR transversal to $\widetilde M$ at $p_0$;

\smallskip
{\rm (ii)}  $H$ is a finite mapping.
\end{Thm}

The implication (ii)$\implies$(i) in Theorem
\ref{t:maincor2} follows immediately from Theorem
\ref{t:maincor}, while (i)$\implies$(ii) is a special case of
Theorem
\ref{t:Hfinite}. (In fact, the implication (i)$\implies$(ii)
holds without the assumption that $M$ is of finite type at $0$ and (ii)$\implies$(i) holds without the assumption that $M$ is essentially finite at $0$. See Section \ref{s:last} and Theorem \ref{t:maincor}.)

In order to formulate another corollary of Theorem
\ref{t:main}, we restrict to the case of real-analytic
manifolds. We recall that the {\it Segre variety} at $p$,
$\Sigma_p$, of a real-analytic generic submanifold
$M$ of $\bC^N$ (given near a point $p_0\in M$ by $\rho=0$) is
the complex manifold (near $p_0$) defined by the equation
\begin{equation}
\rho(Z,\bar p)=0,
\end{equation}
where, by a slight abuse of notation, $\rho(Z,\zeta)$ denotes
also the complexification of $\rho$. Observe that, as is well
known, if $H\colon (\bC^N,p_0)\to (\bC^N,\tilde p_0)$ is a
holomorphic mapping sending $M$ into another real-analytic
generic submanifold $\widetilde M\subset \bC^N$ through $\tilde
p_0$, then $H$ sends $\Sigma_p$ into $\tilde \Sigma_{H(p)}$,
where $\tilde \Sigma_{\tilde p}$ denotes the Segre variety of
$\widetilde M$ at $\tilde p$, for every $p$ sufficiently close to
$p_0$. We write $h$ for the induced holomorphic mapping
$(\Sigma_{p_0}, p_0)\to (\tilde \Sigma_{\tilde
p_0},\tilde p_0)$. We have the following.

\begin{Thm}\Label{t:maincor3} Let $M,\widetilde M\subset\bC^N$ be
real-analytic generic submanifolds of the same dimension
through
$p_0$ and $\tilde p_0$, respectively, with $\widetilde M$ of
finite type at $\tilde p_0$. Suppose that $H\colon (\bC^N,p_0)
\to (\bC^N,\tilde p_0)$ is a holomorphic mapping with
$H(M)\subset \widetilde M$, and let $h$ denote the induced
mapping between Segre varieties $(\Sigma_{p_0}, p_0)\to (\tilde
\Sigma_{\tilde p_0},\tilde p_0)$.  If $h$ is finite, then $H$
is CR transversal to $\widetilde M$ and finite at $p_0$.
Furthermore, the multiplicity of
$H$ equals that of
$h$.
\end{Thm}

Both Theorems \ref{t:maincor2} and \ref{t:maincor3} in the case
where
$M$ and $\widetilde M$ are hypersurfaces follow from results in
\cite{BRgeom}, but the methods of proof given in the latter
work do not generalize to the case of higher codimensional
manifolds.

We should point out that the conclusions of Theorems
\ref{t:maincor},
\ref{t:maincor1}, and \ref{t:maincor2} still hold if the
finite holomorphic mapping $H$ is replaced by a finite CR
mapping
$f$. Recall that a smooth mapping
$f\colon (M,p_0)\to \bC^N$ is called CR if the tangent mapping
$df$ sends the CR bundle
$T^{0,1}M$ into $T^{0,1}\bC^N$. (In particular, the
restriction to $M$ of a holomorphic mapping $H\colon
(\bC^N,p_0)\to \bC^N$ is CR.) A CR mapping $f$ is called {\it
finite} if the associated formal power series mapping  is finite (see Sections
 \ref{s:prel1} and \ref{s:new}). If
$\widetilde M\subset \bC^N$ is another generic submanifold and
$f(M)\subset \widetilde M$, then $f$ is said to be {\it CR
transversal} to  $\widetilde M$ at $p_0$ if
\begin{equation}\Label{e:CRtrans0}
T^{1,0}_{f(p_0)}\widetilde M+T^{0,1}_{f(p_0)}\widetilde M+
df(\bC T_{p_0} M)=\bC T_{f(p_0)} \widetilde M.\end{equation}
One can check (see Section \ref{s:trans}) that if $f=H|_M$,
for some holomorphic mapping $H$, and $f(M)\subset \widetilde M$, then
$f$ is CR transversal to $\widetilde M$ at
$p_0$ if and only if $H$ satisfies (\ref{e:CRtransdef}).

One of the most important tools used in this paper (and a
crucial new ingredient for the proof of Theorem 
\ref{t:maincor}
beyond the hypersurface case) is that of the
iterated Segre mappings as developed by the authors, jointly
with M. S. Baouendi, in the papers
\cite{BER96},
\cite{BER99b}, and \cite{BER03} (see also the book
\cite{BER99a}). These have also played a significant role in
recent work on other aspects of mappings between generic
submanifolds, especially those of codimension greater than one.
We mention here work of Zaitsev \cite{Z1}, \cite{Z2}, S.-Y. Kim
and Zaitsev \cite{KZ01}, and Mir \cite{Mirh}, \cite{Mirg},
\cite{Mi02b}. See also
\cite{BER98}, \cite{BER00}, \cite{BMR02}, \cite{MMZ03}.

The plan of this paper is as follows. In Section
\ref{s:prel1}, we give some preliminary definitions and
results about formal manifolds and mappings, finite type, and
the iterated Segre mappings. One of our main results,
Theorem \ref{t:main}, is then stated and proved in Section
\ref{s:mainresult}. We also give a slightly stronger version
of this theorem in the case of hypersurfaces in Section
\ref{s:hyper}. In the next section, equivalent notions of
transversality are discussed. Results
on nondegeneracy and multiplicity of mappings between
essentially finite manifolds are proved in Section
\ref{s:new}.  In Section
\ref{s:last}, we show how the results in this introduction can
be derived from the theorems proved in the paper; some
remarks and open questions are also given.

\section{Formal manifolds, Finite type, Segre mappings, and
the iterated reflection identity}\Label{s:prel1}

Let
$\bC[[x]]=\bC[[x_1,\ldots,x_k]]$ be the ring of
formal power series in $x=(x_1,\ldots, x_k)$
with complex coefficients. Suppose that
$\rho=(\rho_1,\ldots,
\rho_d)\in\bC[[Z,\zeta]]^d$, where
$Z=(Z_1,\ldots, Z_N)$ and
$\zeta=(\zeta_1,\ldots,\zeta_N)$, satisfies the
reality condition
$
\rho(Z,\zeta)=\bar\rho(\zeta,Z),
$
where $\bar\rho$ is the formal series
obtained from $\rho$ by replacing each
coefficient in the series by its complex
conjugate.
If, in addition, the series
$\rho$ satisfies the condition
$\rho(0)=0$, and
$$
\partial_{(Z,\zeta)}\rho_1(0)\wedge
\ldots\wedge \partial_{(Z,\zeta)}\rho_d(0)\neq 0,
$$
then we say that $\rho$ defines a {\it formal
real submanifold} $M$ of $\bC^N$ through $0$
of codimension $d$ (and dimension $2N-d$). If the defining
series satisfies the stronger condition
$$
\partial_Z\rho_1(0)\wedge
\ldots\wedge \partial_Z\rho_d(0)\neq 0,
$$
we shall say that $M$ is {\it generic}, and refer to
$n:=N-d$ as the {\it CR dimension} of $M$. If $M$ is a
formal real submanifold of codimension $d=1$, then it is
necessarily generic; we shall refer to such an $M$ as a {\it
formal real hypersurface}. These definitions are
motivated by the fact that if in addition the components of
$\rho$ are convergent power series, then the
equations
$\rho(Z,\bar Z)=0$ define a generic real-analytic
submanifold $M$ of $\bC^N$ through $0$.
 Also, if $M$ is a smooth generic
submanifold in $\bC^N$ through $0$, then the
Taylor series at $0$ of a smooth defining
function $\rho(Z,\bar Z)$ of $M$ near $0$, with $\bar Z$ formally
replaced by
$\zeta$, defines a formal
generic submanifold through $0$ (which by a slight abuse of
notation will still be denoted by $M$). These observations
will be used to deduce the results given in the introduction
from the corresponding results for formal real submanifolds.

Let
$H\colon (\bC^N,0)\to (\bC^{N},0)$ be a formal holomorphic (or
simply formal) mapping, i.e.\ $H\in
\bC[[Z_1,\ldots, Z_N]]^{N}$ and each component of
$H(Z)=(H_1(Z),\ldots, H_{N}(Z))$ has no constant term.
If $M$ and $\widetilde M$ are formal real
submanifolds of
$\bC^N$ defined by formal
series
$\rho(Z,\zeta)=(\rho_1(Z,\zeta),\ldots,
\rho_d(Z,\zeta))$ and
$\tilde \rho(Z,\zeta)=(\tilde\rho_1(Z,\zeta),\ldots,
\tilde \rho_{d}(Z,\zeta))$, respectively,
then we say
that the formal mapping $H$, as above, maps
$M$ into $\widetilde M$, denoted $H(M)\subset \widetilde
M$, if
$$
\tilde \rho(H(Z),\bar H(\zeta))=
c(Z,\zeta)\rho(Z,\zeta),
$$
for some $d\times d$
matrix $c(Z,\zeta)$ of formal power series.

If $M\subset
\bC^N$ is a smooth generic submanifold and $f\colon (M,0)\to
(\bC^N,0)$ is a smooth CR mapping, then one may associate to
$f$ a formal mapping $H\colon (\bC^N,0)\to (\bC^N,0)$ as
follows. Let $x$ be  local coordinates on $M$ near $0$ and
$x\mapsto Z(x)$ the local embedding of $M$ into $\bC^N$ near
$0$. Then, there is a unique formal mapping $H$ such that
the Taylor series of $f(x)$ at $0$ equals $H(Z(x))$ as a
power series in $x$ (see e.g.\ \cite{BER99a}, Proposition
1.7.14). Moreover, if $f$ sends $M$ into another smooth
generic submanifold $\widetilde M$, then the induced formal
mapping $H$ sends $M$ into $\widetilde M$ in the sense
described above.

It will be
 convenient to choose {\it normal
coordinates}, $Z=(z,w)$ and
$\zeta=(\chi,\tau)$
with
$z=(z_1,\ldots, z_n)$, $w=(w_1,\ldots,
w_d)$ (so
$n+d=N)$, $\chi=(\chi_1,\ldots,\chi_n)$,
and $\tau=(\tau_1,\ldots,\tau_d)$, in $\bC^N\times\bC^N$ for
$M$ at $0$. By
this we mean
 a formal change of coordinates $Z=Z(z,w)$ with $Z(z,w)$
 a formal invertible mapping $(\bC^N,0)\to (\bC^N,0)$,
and $\zeta=\bar Z(\chi,\tau)$ the corresponding change,
 such that
$$
\rho(Z(z,w),\bar Z(\chi,\tau))=a(z,w,\chi,\tau)
(w-Q(z,\chi,\tau)),
$$
where $a(z,w,\chi,\tau)$ is an invertible $d\times d$
matrix of formal power series, and the components $Q_j$ of the vector valued
$Q\in\bC[[z,\chi,\tau]]^d$ satisfy
\begin{equation}\Label{e:normal}
Q_j(0,\chi,\tau)=
Q_j(z,0,\tau)=\tau_j, \quad j=1,\ldots,d.
\end{equation}
(See Chapter IV.2 of \cite{BER99a}). Here, and in what follows, we use matrix
notation and the convention that the variables $z\in  \bC^n$,
$w\in \bC^d$ are column vectors; in particular, $A^t$ denotes
the transpose of a matrix $A$ and, hence, $Q(z,\chi,\tau)=(Q^1(z,\chi,\tau), \ldots,
Q^d(z,\chi,\tau))^t$ is a column
vector.
For convenience, we
shall simply say that $M$ is defined by the equation
\begin{equation}\Label{e:def} w=Q(z,\bar z,\bar w);
\end{equation}
the reader is referred to \cite{BER99b} or \cite{BER03} for
further definitions and properties related to formal generic
submanifolds and their mappings.

Now let $M$ be a formal
generic submanifold of codimension $d$ through $0$ in
$\bC^{N}$. We let $(z,w)\in
\bC^n\times\bC^d$ be normal coordinates for $M$ at $0$ so
that $M$ is defined at $0$ by
\eqref{e:def}. If $\widetilde M$ is another formal generic
submanifold of codimension
$d$ through $0$ in $\bC^{N}$, and $H\colon
(\bC^N,0)\to (\bC^{N},0)$ a formal mapping sending $M$
to
$\widetilde M$, we write $H=(F,G)$ in the normal coordinates $(\tilde
z,\tilde w)\in
\bC^n\times\bC^d$  for $\widetilde M$. The condition that $H$ sends $(M,0)$ into
$(\widetilde M,0)$ means  that
\begin{equation}\Label{e:basic1} G=\widetilde Q(F,\bar F, \bar G)
\end{equation} and
\begin{equation}\Label{e:basic2} \bar G=\bar{\widetilde Q}(\bar
F,F,G)
\end{equation} where $F=F(z,w)$, $G=G(z,w)$, $\bar F=\bar
F(\chi,\tau)$, $\bar G=\bar G(\chi, \tau)$, whenever $(z,w;
\chi,\tau)$ satisfies
\begin{equation}\Label{e:def2} w=Q(z,\chi,\tau)\end{equation}
or, equivalently, \begin{equation}\Label{e:def3}
\tau=\bar Q(\chi,z,w).
\end{equation} 
 We say that
$H$ is {\it CR transversal} to $\widetilde M$ at $0$ if
\eqref{e:CRtransdef} (with $p_0=\tilde p_0=0$) holds. This is
equivalent to the condition
\begin{equation}\Label{e:detG1}
\det
\displaystyle\frac {\partial G} {\partial w} (0)
\not= 0.
\end{equation}
(See Theorem \ref{t:trans}.)

  For a positive integer $k$, the {\it $k$th Segre mapping} of $M$ at $0$ is the
mapping $v^k: \bC^{kn}\to\bC^N$ defined by
\begin{equation}\Label{e:k-segre}
\bC^{kn}\ni t=(t^1,\ldots,t^k)\mapsto v^k(t):=(t^k,u^k(t^1,\ldots,t^k)), 
\end{equation} 
where $u^k:\bC^{kn}\to \bC^d$ is given inductively by
\begin{equation}\Label{e:uk}u^1(t^1)=0, \ u^k(t^1,\ldots,t^j)=Q(t^k,
t^{k-1},\ov{ u^{k-1}}(t^1,\ldots,t^{k-1})),\ k
\ge 2.\end{equation}  
For example, if
$k>3$ is odd, one has
\begin{equation} u^k(t^1,\ldots, t^k)=Q(t^k,t^{k-1},\bar
Q(t^{k-1},t^{k-2},Q(t^{k-2},t^{k-3},\ldots,\bar
Q(t^2,t^1,0)\ldots).
\end{equation} 
The reader is referred to \cite{BER99b},
\cite{BER03} for the definition and basic properties of the
Segre mappings; we should point out, however, that the
notation in \cite{BER99b} and
\cite{BER03} differ slightly, and in this paper we use the
notation of the latter paper. If $M$
is real-analytic, then the  Segre mappings are
holomorphic and the Segre variety $\Sigma_0$ is parametrized
by the mapping $\bC^n\ni t^1\mapsto v^1(t^1)\in \bC^N$. 
We shall use the convention that $\tilde{}$ over
$Q$,
$v^k$, etc.\ denotes the corresponding objects for $\widetilde M$.

For a formal mapping
 $r: \bC^{l}\to \bC^m$, we denote by $\Rk r$  the
generic rank of
$r,$ i.e.\ the rank of the $m \times l$ matrix
$\partial r/\partial (s^1,\ldots, s^l)$ over the field of
fractions of the ring of formal power series $\bC[[
s^1,\ldots,s^l ]]$.  Hence $\Rk r = j$ if and only if the
matrix $\partial r/\partial (s^1,\ldots, s^l)$ has a $j\times
j$ minor that is nonvanishing as a power series, but no such
$(j+1)\times (j+1)$ minor. 

We shall write
$\bC T_0 M$, $T^{1,0}_0 M$, and $T^{0,1}_0M$ for the vector
spaces of all complex tangent vectors, all $(1,0)$ tangent
vectors, and all $(0,1)$ tangent vectors, respectively, at
$0$, with analogous notation for $\widetilde M$. Recall that $M$ is
said to be of {\it finite type} at
$0$ (in the sense of Kohn, and Bloom--Graham) if the Lie algebra
$\frak g_M$ generated by all formal $(1,0)$ and $(0,1)$
vector fields tangent to $M$ satisfies $\frak
g_M(0)=\bC T_0M$. We shall need the following result
from the joint work of Baouendi and the authors. 
\begin{Thm}[\cite{BER99b}, Theorem 3.1.9] \Label{t:basic} Let
$M\subset
\bC^N$ be a formal generic submanifold of codimension $d$
through
$0$ and
 $v^k(t) = (t^k,u^k(t)):\bC^{k(N-d)}\to \bC^N$ the $k$th
 Segre mapping of $M$, and $\Rk$ the generic rank.
  Then the following are equivalent.

{\rm (i)} $M$ is of finite type at $0$.
\smallskip

{\rm (ii)} $\Rk v^k = N$ for $k\geq d+1$.
\smallskip

{\rm (iii)} $\Rk u^k = d$ for $k\geq d+1$.
\end{Thm} 

Note that (cf.\ Proposition 3.1.6 in
\cite{BER99b}) that
\eqref{e:basic1} and \eqref{e:basic2} hold with
\begin{equation}\Label{e:v1}
(z,w;\chi,\tau)=(v^{l}(t^1,\ldots,t^{l});\ov
{v^{l-1}}(t^1,\ldots, t^{l-1}))\end{equation} or
\begin{equation}\Label{e:v2}
(z,w;\chi,\tau)=(v^{l-1}(t^1,\ldots,t^{l-1});\ov
{v^{l}}(t^1,\ldots, t^{l})), 
\end{equation}
for any $l \ge 2$. 
Let us fix $k$ and substitute \eqref{e:v1} with $l=k$ in
\eqref{e:basic1}. We then obtain the formal power series identity

\begin{equation}
G\circ v^k=\widetilde Q(F\circ v^k,\overline{ F\circ  v^{k-1}},\overline {
G\circ  v^{k-1}}).
\end{equation}
Next, we substitute \eqref{e:v2} with $l=k-1$ in
\eqref{e:basic2} to obtain
\begin{equation}
\overline{ G\circ  v^{k-1}}= \overline{\widetilde Q}(\overline{ F\circ
 v^{k-1}}, F\circ
 v^{k-2},G\circ v^{k-2}),
\end{equation}
and then substitute for $\overline{ G\circ  v^{k-1}}$ in the
previous equation. Continuing this process inductively (and
using the well known fact that, in normal coordinates, $G\circ
v^1(t^1)=G(t^1,0)=0$), we deduce that, for any $k\geq 1$ odd,
\begin{equation}\Label{e:oddid} G\circ v^{k}=\tilde
u^k(F\circ v^1,\ov{F\circ v^2},F\circ v^3,\ldots, F\circ v^k)
\end{equation} and, for $k\geq 2$ even,
\begin{equation}\Label{e:evenid} G\circ v^{k}=\tilde
u^k(\ov{F\circ v^1},F\circ v^2,\ov{F\circ v^3},\ldots, F\circ
v^k),
\end{equation} where $F\circ v^j=(F\circ
v^j)(t^1,\ldots,t^j)$, $G\circ v^j=(G\circ
v^j)(t^1,\ldots,t^j)$ and $\ov{F\circ v^j}=(\ov{F\circ
v^j})(t^1,\ldots,t^j)$; recall also that $\tilde v^j(\tilde
t^1,\ldots,\tilde t^j)=(\tilde t^j,\tilde u^j(\tilde
t^1,\ldots,\tilde t^j))$ denotes the Segre mapping of the
target manifold $\widetilde M$. It follows from \eqref{e:oddid}
that (for $k\geq 1$ odd)
\begin{equation}\Label{e:oddid2} H\circ v^{k}=\tilde
v^k(F\circ v^1,\ov{F\circ v^2},F\circ v^3,\ldots, F\circ v^k)
\end{equation} and from \eqref{e:evenid} that (for $k\geq 2$
even) 
\begin{equation}\Label{e:evenid2} H\circ v^{k}=\tilde
v^k(\ov{F\circ v^1},F\circ v^2,\ov{F\circ v^3}\ldots, F\circ
v^k).
\end{equation}
We shall refer to the formal power series identities \eqref{e:oddid2} and
\eqref{e:evenid2} (and also
\eqref{e:oddid} and \eqref{e:evenid}) as {\it iterated
reflection identities}.

\begin{Rem} {\rm Note that if $M$ is of finite type, the mapping $H$ is completely determined by the
component $F$. Indeed, by Theorem \ref{t:basic}, if $k \ge
d+1$, then $v^k$ is of generic rank $N$. Hence  for $k$ large,
$H$ is determined by the right hand side of 
\eqref{e:oddid2} (or
 \eqref{e:evenid2}), which depends only on $F$.}\end{Rem}

Recall that the mapping $H=(F,G)$ is called {\it not totally
degenerate} at $0$ if
\begin{equation}\Label{e:Hntd}
\det \frac{\partial F}{\partial z}\circ v^1\not\equiv 0.
\end{equation}  If $H$ is holomorphic, then
$H$ is not totally degenerate at $0$ if and only if the Jacobian
determinant of the induced mapping
$h\colon (\Sigma_0,0)\to (\tilde \Sigma_0,0)$ is not
identically zero. We shall use the notation $\Jac H$ for the
Jacobian determinant of $H$, i.e.\ $\Jac H=\det(\partial
H/\partial Z)$ where $Z$ is a coordinate on $\bC^N$. We have
the following result.

\begin{Pro}\Label{p:ft} Let $M,\widetilde M$ be formal generic
submanifolds of the same dimension through $0\in \bC^N$,
and
$H\colon (\bC^N,0)\to (\bC^N,0)$ a formal holomorphic mapping
with $H(M)\subset\widetilde M$. Then the following
hold:
\smallskip

{\rm (a)} If $\Jac H\not \equiv 0$  and $M$ is of finite
type at
$0$, then $\widetilde M$ is of finite type at $0$.
\smallskip

{\rm (b)} If $H$ is not totally degenerate at $0$ and $\widetilde M$
is of finite type at $0$, then $M$ is of finite type at $0$ and
$\Jac(H)\not\equiv 0$.
\end{Pro}

\begin{proof} To prove (a), assume that $M$ is of finite
type at
$0$. Hence  $\Rk v^k=N$ for $k\geq d+1$ by Theorem
\ref{t:basic}. If $\Jac H\not\equiv 0$, i.e.\ $\Rk H =N$,
then the generic rank of the left hand side of
\eqref{e:oddid2} (and
\eqref{e:evenid2}), for $k\geq d+1$, is $N$. It follows
that 
$\Rk \tilde v^k = N$, $k\geq d+1$,  and hence
$\widetilde M$ is of finite type at
$0$, again by Theorem \ref{t:basic}.

To prove (b), let us choose $k$ odd, and define $\Psi\colon
(\bC^{kn},0)\to (\bC^{kn},0)$ by
\begin{equation}
\Psi(t^1,\ldots, t^k)=(F\circ v^1(t^1), \overline{F\circ
v^2}(t^1,t^2),\ldots, F\circ v^k(t^1,\ldots, t^k)).
\end{equation} Observe that the $kn\times kn$ matrix
$\partial \Psi/\partial t$, with $t:=(t^1,\ldots, t^{k})$, is
lower block triangular. That is, there are $n\times n$
 diagonal block matrices
$D_1,\ldots, D_k$ such that the $D_j$ are on the diagonal  
of $\partial \Psi/\partial t$ with all entries zero above the diagonal
blocks.  The blocks $D_j(t)$ are given by 
\begin{equation}\begin{aligned} D_1(t) &={\frac{\partial
F}{\partial z}\circ v^1}(t^1)\\ D_2(t) &= \ov{\frac{\partial
F}{\partial z}\circ v^2}(t^1,t^2)+\ov{\frac{\partial
F}{\partial w}\circ v^2}(t^1,t^2)
\ov{ \frac{\partial
Q}{\partial z}}(t^2,v^1(t^1))\\ & \vdots\\
 D_k(t) &=
{\frac{\partial F}{\partial z}\circ v^m}(t^1,\ldots,
t^m)+{\frac{\partial F}{\partial w}\circ v^m}(t^1,\ldots,
t^m)
 {\frac{\partial  Q}{\partial
z}}(t^m,\ov{v^{m-1}}(t^1,\ldots, t^{m-1})).\\
\end{aligned}
\end{equation} 
 Since $(\partial Q/\partial z)(z,0,0)=0$, it follows that
if 
$t^1=\ldots=t^{j-1}=0$, then $$D_j(t)=(\partial F/\partial
z\circ v^1)(t^j),\ j\  \text{\rm odd},\ \  
D_j(t)=(\ov{\partial F/\partial z\circ v^1})(t^j),\  j\
\text{\rm even}.$$
 We conclude, since $H$ is assumed not totally
degenerate (see \eqref{e:Hntd}), that $\partial \Psi/\partial
t$ is invertible over the field of fractions of formal power
series in $t$. If
$\widetilde M$ is assumed to be of finite type at 0, then
$\Rk\tilde v^k=N$  for $k\geq d+1$ by Theorem
\ref{t:basic}, so that the right hand side of
\eqref{e:oddid2} (which is equal to $\tilde v^k\circ\Psi$ in this
notation) has generic rank $N$. Hence, the left hand
side of \eqref{e:oddid2} must also have generic rank
$N$. It follows that both
$\Rk H=N$ (i.e. $\Jac H \not\equiv 0$) and
$\Rk v^k=N$. The latter implies that $M$ is of finite type at 0, again by
Theorem \ref{t:basic}. This completes the proof.
\end{proof}

We would like to point out that the condition ``$\Jac
H\not\equiv0$" in statement (a) of Proposition \ref{p:ft}
cannot be replaced by the condition ``$H$ is not
totally degenerate", nor can the latter condition in statement
(b) be replaced by the former as is illustrated by the
following two examples.

\begin{Ex} {\rm Let $M\subset \bC^2$ be the real hypersurface
given by
$
\im w=|z|^2,
$
and $\widetilde M\subset \bC^2$ given by
$
\im \tilde w=\re \tilde w|\tilde z|^2.
$
The mapping $H\colon (\bC^2,0)\to (\bC^2,0)$ given by
$
H(z,w)=(z,0)
$
is not totally degenerate and sends $M$ into $\widetilde M$.
Note that $M$ is of
finite type at $0$ but $\widetilde M$ is of infinite type at $0$.}
\end{Ex}

\begin{Ex} {\rm Let $M\subset \bC^2$ be the real hypersurface
given (in implicit form) by
$
\im w=|zw|^2,
$
and $\widetilde M\subset \bC^2$ given by
$
\im \tilde w=|\tilde z|^2.
$
The mapping $H\colon (\bC^2,0)\to (\bC^2,0)$ given by
$
H(z,w)=(zw,w)
$
satisfies $\Jac H\not\equiv 0$ and sends $M$ into $\widetilde M$.
Note that $\widetilde M$ is of
finite type at $0$ but $M$ is of infinite type at $0$.}
\end{Ex}

We conclude this section by giving two definitions for a
formal holomorphic mapping $H\colon (\bC^N,0)\to (\bC^N,0)$
and a lemma describing how these notions are related. We say
that $H$ is {\it finite} if
$$
\dim_\bC \bC\dbl Z\dbr/I(H(Z))<\infty,
$$
where $I(H(Z))$ denotes the ideal generated by the components
of the mapping $H(Z)=(H_1(Z),\ldots, H_N(Z))$. In particular,
if $H$ is finite, then there is a positive integer $k$ such
that any holomorphic mapping $K\colon (\bC^N,0)\to
(\bC^N,0)$ with $H(Z)=K(Z)+O(|Z|^{k+1})$ is a finite
holomorphic mapping (see e.g.\ Proposition 5.1.8 in
\cite{BER99a}). 

For our last definition, we assume again that $H$ sends $M$
into $\widetilde M$ and write $H=(F,G)$ in
the normal coordinates $(\tilde z,\tilde w)$ for
$\widetilde M$. We say that $H$ is {\it Segre finite} if the
formal mapping $F\circ v^1\colon (\bC^n,0)\to (\bC^n,0)$ is
finite. We have the following relations between the notions
of nondegeneracy defined above.

\begin{Lem}\Label{l:wellknown0} Let $M,\widetilde M$ be
formal generic submanifolds of the same CR dimension through
$0\in
\bC^N$, and
$H\colon (\bC^N,0)\to (\bC^N,0)$ a formal holomorphic mapping
sending $M$ into $\widetilde M$. The following hold.
\medskip

{\rm (i)}
$H$ is finite $\implies$ $H$ is Segre finite $\implies$ $H$
is not totally degenerate.
\smallskip

{\rm (ii)} If $H$ is Segre finite and CR transversal, then
$H$ is finite.
\end{Lem}

The lemma is a consequence of Lemma \ref{l:wellknown} below.

\section{CR transversality of finite
mappings}
\Label{s:mainresult}

In this section, we prove, in the context of formal manifolds and
mappings, one of the main results of this paper. Our theorem is the following. 

\begin{Thm}\Label{t:main} Let $M,\widetilde M$ be formal 
generic submanifolds of the same dimension through $0\in
\bC^N$ with $\widetilde M$
of finite type at
$0$, and
$H\colon (\bC^N,0)\to (\bC^N,0)$ a formal holomorphic mapping
with $H(M)\subset\widetilde M$.  If $H$
is a finite map, then $H$ is CR
transversal to $\widetilde M$ at $0$. More generally, 
if $H$ is Segre finite at $0$, then $H$ is a finite mapping
that is CR transversal to $\widetilde M$ at $0$. 
\end{Thm}

In the special case that $M$ and $\widetilde M$ are hypersurfaces, Theorem
\ref{t:main} is proved in \cite{BRgeom}. Our approach here,
using the iterated Segre mappings, is different from that of
\cite{BRgeom}. In fact, we obtain a slightly stronger
result in the hypersurface case than that proved in
\cite{BRgeom} (see Theorem \ref{t:hyper} below).

\begin{Rem} {\rm Without further
assumptions on $M$, there may exist a holomorphic mapping $H$ sending $M$ into $\widetilde M$ such that $H$ is CR transversal to $\widetilde M$ at $0$ but not Segre finite. However, if $M$ is essentially
finite at $0$, then every CR transversal holomorphic mapping is necessarily finite and, hence, also Segre finite (see Theorem  \ref{t:Hfinite} and Lemma \ref{l:wellknown}).  }
\end{Rem}

To prove Theorem \ref{t:main}, it suffices, by Lemma
\ref{l:wellknown0}, to prove that if
$H$ is Segre finite, then it is CR transversal. For the
proof we begin with some lemmas. In these lemmas,
$M$ and $\widetilde M$ will be formal generic submanifolds
of the same dimension through
$0\in \bC^N$, and
$H\colon (\bC^N,0)\to (\bC^N,0)$ a formal holomorphic mapping
sending $M$ into $\widetilde M$.
 We choose normal coordinates
for $M$ and
$\widetilde M$, respectively, and
write $H=(F,G)$ as in Section \ref{s:prel1}. We shall use the
notation $v^k$ for the Segre mappings introduced in
\eqref{e:k-segre} and denote the corresponding Segre mappings
for $\widetilde M$ by $\tilde v^k$. We fix an integer
$m\geq d+1$ 
and for $t^i, \tilde
t^i \in \bC^n$, $1\le i \le 2m$, we write
$t=(t^1,\ldots,t^{2m})$, $t'=(t^1,\ldots,t^{2m-1})$, $\tilde
t=(\tilde t^1,\ldots,\tilde t^{2m})$,$\tilde
t'=(\tilde t^1,\ldots,\tilde t^{2m-1})$.  We introduce the formal mapping $\Phi\colon
(\bC^{(2m-1)n},0)
\to (\bC^{(2m-1)n},0)$ given by
\begin{equation}\Label{e:Phi}
\Phi(t^1,\ldots, t^{2m-1}):=(\ov{F\circ v^1}(t^1),F\circ
v^2(t^1,t^2),\ldots, \ov{F\circ v^{2m-1}}(t^1,\ldots,
t^{2m-1})).
\end{equation}
We also define a
subspace $W \subset \bC^{2mn}$ as follows.
\begin{equation}\Label{e:vars} W:= \{ t \in \bC^{2mn}:
t^{2m}=0,\ t^{2m-1}=t^1,\ t^{2m-2}=t^2,\ldots,\
t^{m+1}=t^{m-1}\}.\end{equation}
The subspace $W\subset \bC^{2mn}$ may also be 
identified with a subspace $W'\subset\bC^{(2m-1)n}$,
where
\begin{equation}\Label{e:varsp} W' =\{t' \in \bC^{(2m-1)n}: 
(t',0) \in
W\}.\end{equation}
 One can check, by using \eqref{e:uk} with $k = 2m$ and the
identities
\begin{equation} Q(z,\chi,\bar Q(\chi,z,w))=w,\ \bar
Q(\chi,z,Q(z,\chi,\tau))=\tau,
\end{equation} as well as the corresponding ones for $Q$ replaced by
$\widetilde Q$ (cf.\ also Lemma 4.1.3 in
\cite{BER99b}), that one has
\begin{equation}\Label{e:tildevars1} 
\begin{aligned}
v^{2m}|_W = 
\tilde v^{2m}|_{W} = 0,\ & \ u^{2m-1}|_{W'} = 
\tilde u^{2m-1}|_{W'}=0\\  \Phi(W') \subset & 
W'.
\end{aligned}
\end{equation}
\begin{Rem}\Label{r:rank}{\rm The significance of the 
subspace $W'$
given by \eqref{e:varsp} is expressed in the following result proved
in 
\cite{BER99b}, Theorem 3.1.9 and Lemma 4.1.3 (see also
\cite{BER03}, Theorem 2.4):}
\medskip

\noindent {\it The generic submanifold $M$ (resp.\ $\widetilde M$) is of finite type at
$0$ if and only if the $d\times (2m-1)n$-matrix
$\displaystyle
\frac{\partial u^{2m}}{\partial t'}(t',0)$ (resp.\
$\displaystyle \frac{\partial
\tilde u^{2m}}{\partial \tilde t'}(t',0)$)  has a $d\times
d$ minor that is nonvanishing on $W'$.}
\end{Rem}
\medskip

 We may now state our
first lemma, which gives the main identity
used in the proof of Theorem \ref{t:main}.

\begin{Lem}\Label{l:mainid}   Let $W\subset \bC^{2mn}$ be
given by \eqref{e:vars}. Then the identity
\begin{equation}\Label{e:diffid}
\frac{\partial G}{\partial w}(0)\frac{\partial
u^{2m}}{\partial t'}(t) =\frac{\partial \tilde
u^{2m}}{\partial
\tilde t'}(\Phi(t'),0)\frac{\partial \Phi}{\partial
t'}(t')+\frac{\partial \widetilde Q}{\partial \tilde
z}(0,\ov{F\circ v^1}(t^1),0)\frac{\partial F}{\partial
w}(0)\frac{\partial u^{2m}}{\partial t'}(t),
\end{equation}
holds for all $t = (t',0) \in W$.
\end{Lem}

\begin{proof}  We observe that \eqref{e:evenid} with $k=2m$ can
be written as
\begin{equation}\Label{e:help}
G\circ v^{2m}(t)=\tilde u^{2m}(\Phi(t'),F\circ v^{2m}(t)),
\end{equation}
where $\Phi(t')$ is given by \eqref{e:Phi}.
We differentiate \eqref{e:help} in the first set
of
$(2m-1)n$ variables $t'=(t^1,\ldots, t^{2m-1})$, and then
restrict to the subspace $W$.  We
obtain the desired identity \eqref{e:diffid} by using 
\eqref{e:tildevars1}, as well as the chain rule.
\end{proof}

\begin{Lem}\Label{l:A} Assume that $\widetilde M$ is of finite
type at
$0$ and that $H$ is not totally degenerate. If there is $V\in \bC^d
\setminus \{0\}$  such that
$
V^t\displaystyle\frac{\partial G}{\partial w}(0)=0,
$
then
\begin{equation}\Label{e:contra0} 
V^t\frac{\partial \widetilde Q}{\partial \tilde
z}(0,\ov{F}(t^1,0),0)\frac{\partial F}{\partial
w}(0)\not\equiv 0.
\end{equation}
In particular, we have
\begin{equation}\Label{e:contra} V^t\frac{\partial \widetilde Q}{\partial \tilde z}(0,\ov{F}(t^1,0),0)\not\equiv 0.
\end{equation}
\end{Lem}

\begin{proof} We shall show, under the hypotheses of the
lemma, that the $d\times (2m-1)n$ matrix
\begin{equation}\Label{e:prod}
\frac{\partial \tilde u^{2m}}{\partial \tilde
t'}(\Phi(t'),0)\frac{\partial \Phi}{\partial t'}(t')
\end{equation} 
has a $d\times d$ minor
that does not vanish identically on $W'$ (given by \eqref{e:varsp}). The
conclusion of the lemma will then easily follow from the
identity
\eqref{e:diffid} by multiplying both sides to the left by $V^t$.

A direct calculation shows
that the $(2m-1)n\times (2m-1)n$ matrix
$\displaystyle\frac{\partial
\Phi}{\partial t'}(t')$ is a lower triangular block matrix,
with
 $n\times n$-matrices $D_1,D_2,\ldots, D_{2m-1}$ on
the diagonal and zeroes above.  More precisely, for
$t'=(t^1,\ldots, t^{2m-1}) \in \bC^{(2m-1)n}$ we have
\begin {equation} D_k(t') = \ov{\frac{\partial  (F\circ
v^k)}{\partial t^k}}(t^1,\ldots,t^k)\ \text{\rm for $k$ odd,}\quad 
D_k(t') = {\frac{\partial  (F\circ v^k)}{\partial t^k}}(t^1,\ldots,t^k)\
\text{\rm for $k$ even,}
\end{equation}
$k =
1,\ldots, 2m-1 $.  Hence we obtain
\begin{equation}\begin{aligned} D_1(t')
&={\ov{\frac{\partial F}{\partial z}\circ v^1}}(t^1)\\
D_2(t') &= {\frac{\partial F}{\partial z}\circ
v^2}(t^1,t^2)+{\frac{\partial F}{\partial w}\circ
v^2}(t^1,t^2) { \frac{\partial
Q}{\partial z}}(t^2,\ov {v^1}(t^1))\\ & \vdots\\
 D_{2m-1}(t') &=
\ov{{\frac{\partial F}{\partial z}\circ
v^{2m-1}}}(t')+\ov{{\frac{\partial F}{\partial w}\circ
v^{2m-1}}}(t')
 \ov{{\frac{\partial  Q}{\partial
z}}}(t^{2m-1},v  ^{2m-2} (t^1,\ldots, t^{2m-2})).\\
\end{aligned}
\end{equation} For each $j$, $1 < j\leq
2m-1$, let $W'_j\subset W'$ be defined by
\begin{equation} W_j'= \big\{(t^1,\ldots, t^{2m-1})\in
\bC^{(2m-1)n}: t^k= 0\ \text{\rm if}\ k\not= j \ \text{\rm
or}\ k\not=2m-j, \ \text{\rm and}\ t^j = t^{2m-j}\}.
\end{equation}
Since $Q(z,0,0)\equiv Q(0,\chi,0) \equiv 0$, one can
check that for
$t'=(t^1,\ldots, t^{2m-1})\in W_j'$,
$$v^j(t^1,\ldots,t^j)\equiv v^1(t^j),\  v^{j-1}(t^1,\ldots,t^{j-1})
\equiv 0.$$
Hence, for $t'\in W'_j$, we have $D_j(t')=(\partial
F/\partial z)
\circ v^1(t^j)$, or its complex conjugate, depending on the
parity of $j$. Since
$H$ is not totally degenerate, i.e.\
$\bC^{n}\ni t^1\mapsto\displaystyle\det\frac{\partial F}{\partial z}\circ v^1(t^1) 
\not\equiv 0,$ it follows that 
$\displaystyle\det\frac{\partial
\Phi}{\partial t'}(t')|_{W'}\not\equiv 0$.  
 By a similar calculation, we can also check that 
$\Phi|_{W'}: W' \to W'$ (see
\eqref{e:tildevars1}) is of full generic rank. Finally,
since
$\widetilde M$ is assumed to be of finite type at $0$,
  the $d\times (2m-1)n$ matrix $\displaystyle\frac {\partial
\tilde u^{2m}}{\partial \tilde t'}(\tilde t',0)|_{W'}$ has a nonvanishing $d\times
d$ minor (see Remark \ref{r:rank}). It follows that  
the
$d\times (2m-1)n$ matrix
 $\displaystyle\frac {\partial
\tilde u^{2m}}{\partial  \tilde t'}(\Phi(\tilde t'),0)|_{W'}$ also  has a
nonvanishing
$d\times d$ minor.    Hence the product in
\eqref{e:prod} has a nonvanishing
$d\times d$ minor on $W'$. This completes
the proof of Lemma \ref{l:A}. \end{proof}

We may now reduce the proof of Theorem
\ref{t:main} to our third lemma (Lemma \ref{l:B} below).
\begin{proof}[Proof of Theorem $\ref{t:main}$] We must prove
that $\det \partial G/\partial w(0)\neq 0$ (see Theorem
\ref{t:trans}). Let us assume, in order to reach a
contradiction, that
$\partial G/\partial w(0)$ is not invertible. Since, as is
well known, $\partial G/\partial w(0)$ is real, it follows
that there is nonzero $V\in \bR^d$ such that
\begin{equation}\Label{e:ass} V^t\frac{\partial G}{\partial
w}(0)=0.
\end{equation} We then conclude, in view of Lemma \ref{l:A}
and the fact that Segre finite maps are not totally
degenerate (see Lemma \ref{l:wellknown}), that
\begin{equation}\Label{e:contra2} V^t\frac{\partial \widetilde Q}{\partial \tilde z}(0,\ov{F}(t^1,0),0)\not\equiv 0.
\end{equation}
We shall show that \eqref{e:contra2} is impossible under the assumptions of the
theorem. To this end, we go back to \eqref{e:evenid} with $k=4$ and complex
conjugate this identity to obtain
$$
\ov{G\circ v^4}=\ov{\tilde u^4}(F\circ v^1,\ov{F\circ
v^2},F\circ v^3,\ov{F\circ v^4}).
$$ 
Proceeding as in Lemma \ref{l:mainid},  we differentiate this
identity with respect to $t^3$ and then set $t^4=0$, $t^3=t^1$
(i.e. restrict to $W$ in the case $m=2$). We obtain
\begin{multline}\Label{e:newid}
\ov{\frac{\partial G}{\partial w}}(0)\frac{\partial
Q}{\partial z}(t^1,\ov{v^2})=\frac{\partial \ov{\widetilde Q}}{\partial \tilde
\chi } (0,F\circ v^1,0)\ov{\frac{\partial F}{\partial
w}}(0)\frac{\partial Q}{\partial z}(t^1,\ov{v^2})+\\
\frac{\partial \widetilde Q}{\partial \tilde z}(F\circ
v^1,\ov{F\circ v^2},\ov{\widetilde Q}(\ov{F\circ v^2},F\circ
v^1,0))\left(\frac{\partial F}{\partial z}\circ
v^1+\frac{\partial F}{\partial w}\circ v^1\, \frac{\partial
Q}{\partial z}(t^1,\ov{v^2})\right ),
\end{multline} where $v^1=v^1(t^1)$ and
$\ov{v^2}=\ov{v^2}(t^1,t^2)$. (Note that \eqref{e:newid} is the
complex conjugate of one of the components of \eqref{e:diffid}
with
$m=2$.) Also, observe that
$$\ov{\widetilde Q}(\ov{F\circ v^2}(t^1,t^2),F\circ
v^1(t^1),0)=\ov{ G}(t^2,\ov
Q(t^2,t^1,0)),$$
by (the complex conjugate of) \eqref{e:evenid} with $k=2$.
For convenience, we shall
replace the variables $(t^1,t^2)$ in
\eqref{e:newid} by $(z,\chi)$.
If we multiply \eqref{e:newid}
from the left by $V^t$ and use the assumption \eqref{e:ass},
we then obtain from \eqref{e:newid}
\begin{multline}\Label{e:LTid} V^t\frac{\partial \ov{\widetilde Q}}{\partial \tilde
\chi}(0,F(z,0),0)\ov{\frac{\partial F}{\partial
w}}(0)\frac{\partial Q}{\partial z}(z,\chi,\bar
Q(\chi,z,0))=\\- V^t\frac{\partial \widetilde Q}{\partial
\tilde z}(F(z,0),\bar F(\chi,\ov Q(\chi,z,0)), \ov
G(\chi,\ov Q(\chi,z,0)))\times \\
\left(\frac{\partial F}{\partial z}(z,0)+\frac{\partial
F}{\partial w}(z,0)\frac{\partial Q}{\partial z}(z,\chi,\ov
Q(\chi,z,0))\right ).
\end{multline}
The contradiction needed to complete the  proof of
Theorem
\ref{t:main} is now obtained by the following lemma, in
view of Lemma
\ref{l:A}.
\end{proof}

\begin{Lem}\Label{l:B} Assume that $H$ is Segre finite. If
$V\in \bC^d$ is a nonzero vector for which \eqref{e:LTid}
holds, then
\begin{equation}\Label{e:contra3} V^t\frac{\partial \widetilde Q}{\partial \tilde z}(0,\ov{F}(\chi,0),0)\equiv 0.
\end{equation}
\end{Lem}
\begin{proof} Suppose, in order to reach a contradiction,
\eqref{e:contra3} is false. Then let
$\alpha_0$ be a multi-index of minimal length such that
there is a nonzero vector $U\in \bC^n$ such that 

\begin{equation}\Label{e:contra4}
\bigg(\frac{\partial}{\partial\chi}\bigg)^{\alpha_0}
\bigg(V^t\frac{\partial
\widetilde Q}{\partial \tilde z}(0,\ov{F}(\chi,0),0)\bigg)=
 U^t + O(|\chi|).
\end{equation}

To complete the proof of Lemma \ref{l:B}, we shall need the
following general fact about finite formal maps.

\begin{Lem}\Label{l:finite} If $f:(\bC^n,0) \to (\bC^n,0)$ is
a finite formal map and $\beta \in \bC^n \setminus \{0\}$ is
any nonzero vector, then there is a formal curve
$\gamma: (\bC,0) \to (\bC^n,0)$ and a positive integer $p$
such that 
\begin{equation}\Label{e:fgamma} f(\gamma(t)) = t^p\beta +
O(t^{p+1}).
\end{equation} 
\end{Lem}
\begin{proof} Let $L$ be the
complex line $\{\tilde z\in \bC^n\colon \tilde z=s\beta\ \text{\rm with }
s\in \bC\}$ and $l_1(\tilde z),\ldots,
l_{n-1}(\tilde z)$ linear functions generating the ideal of
$L$. Let $\tilde J$ the ideal in $\bC[[\tilde z]]$
generated by $l_1(\tilde z),\ldots,
l_{n-1}(\tilde z)$ and $J$ the ideal generated by
$(l_i\circ f)(z)$, $i=1,\ldots,n-1$. Since $J$ is generated by $n-1$ elements,
$J$ cannot have finite codimension in $\bC[[z]]$. It
follows (see e.g.\ Lemma 3.32 of \cite{BER00}) that there
exists a nontrivial formal curve $\gamma\colon (\bC,0)\to
(\bC^n,0)$ such that $g\circ\gamma\equiv 0$ for every $g\in
J$. In particular, we have $l_i\circ f\circ \gamma\equiv 0$
for $i=1,\ldots,n-1$. Moreover, since $f$ is a finite formal
mapping, it also follows Lemma 3.32 of \cite{BER00} that
$f\circ\gamma\not\equiv 0$. 
Expanding $f(\gamma(t))=t^p\beta'+O(t^{p+1})$, with
$\beta'\neq 0$, we obtain
$0\equiv l_i(\beta')t^p+O(t^{p+1})$ for $i=1,\ldots, n-1$
and, hence, we conclude that $\beta'\in L$. This completes
the proof.
\end{proof}

We may now complete the proof of Lemma \ref{l:B} and, hence, 
that of Theorem \ref{t:main}. Let $h(t):=F(\gamma(t),0)$,
where
$\gamma\colon (\bC,0)\to (\bC^n,0)$ is the complex analytic
curve obtained by applying Lemma \ref{l:finite} with
$f(z):=F(z,0)$ and $\beta=\bar U$. Then,
\begin{equation}\Label{e:direction}F_z(\gamma(t),0)\gamma'(t)=h'(t)=pt^{p-1}\bar
U+O(t^p).
\end{equation}
 We set
$z=\gamma(t)$ in \eqref{e:LTid}, multiply both sides of
the resulting identity by $\gamma'(t)$,  apply
$(\partial/\partial\chi)^{\alpha_0}$, where
$\alpha_0$ is as in \eqref{e:contra4}, and set $\chi=0$. We
obtain
\begin{multline}\Label{e:LTid2}
\bigg(\frac{\partial}{\partial
\chi}\bigg)^{\alpha_0}\bigg(V^t\frac{\partial
\ov{\widetilde Q}}{\partial \tilde
\chi}(0,h(t),0)\ov{\frac{\partial F}{\partial
w}}(0)\frac{\partial Q}{\partial z}(\gamma(t),\chi,\bar
Q(\chi,\gamma(t),0))\gamma'(t)\bigg)\bigg|_{\chi=0}=\\-
\bigg(\frac{\partial}{\partial
\chi}\bigg)^{\alpha_0}\bigg(V^t\frac{\partial
\widetilde Q}{\partial
\tilde z}(h(t),\bar F(\chi,\ov Q(\chi,\gamma(t),0)), \ov
G(\chi,\ov Q(\chi,\gamma(t),0)))\times \\
\big(h'(t)+\frac{\partial
F}{\partial w}(\gamma(t),0)\frac{\partial Q}{\partial
z}(\gamma(t),\chi,\ov Q(\chi,\gamma(t),0))\gamma'(t)\big
)\bigg)\bigg|_{\chi=0}.
\end{multline} 
 Since $
V^t(\partial \ov{\widetilde Q}/\partial \tilde
\chi)(0,0,0)=0
$
and $h(t)=O(t^p)$,
we conclude that the left hand side of \eqref{e:LTid2} is
$O(t^p)$. Also, observe that by \eqref{e:contra4} we have
\begin{equation}
\bigg(\frac{\partial}{\partial
\chi}\bigg)^{\alpha_0}\bigg(V^t\frac{\partial
\widetilde Q}{\partial
\tilde z}(h(t),\bar F(\chi,\ov Q(\chi,\gamma(t),0)), \ov
G(\chi,\ov Q(\chi,\gamma(t),0)))\bigg)\bigg|_{t=0}=U^t +
O(|\chi|)
\end{equation}
and 
\begin{equation}
\bigg(\frac{\partial}{\partial
\chi}\bigg)^{\alpha}\bigg(V^t\frac{\partial
\widetilde Q}{\partial
\tilde z}(h(t),\bar F(\chi,\ov Q(\chi,\gamma(t),0)), \ov
G(\chi,\ov Q(\chi,\gamma(t),0)))\bigg)\bigg|_{t=0}=0,
\end{equation}
for all $|\alpha|<|\alpha_0|$. Hence, by 
\eqref{e:direction} and \eqref{e:normal},
the right hand side of \eqref{e:LTid2} is
$p|U|^2t^{p-1}+O(t^p)$. This is a contradiction
since $U\neq 0$. This completes the proof of Lemma \ref{l:B}.
\end{proof}

\section{CR transversality in the hypersurface
case}\Label{s:hyper}

We shall in this section give a slight improvement of Theorem
\ref{t:main} in the case where $M$ and $\widetilde M$ are
hypersurfaces (i.e.\ $d=1$). Recall that if a formal mapping $H$ sending $M$ into $\widetilde M$
 is Segre finite, then $H$  is also not totally degenerate (see
Lemma \ref{l:wellknown0}).  The converse is not true in
general. 

\begin{Thm}\Label{t:hyper} Let $M,\widetilde M$ be formal
real hypersurfaces through $0\in
\bC^N$, and
$H\colon (\bC^N,0)\to (\bC^N,0)$ a formal holomorphic mapping
with $H(M)\subset\widetilde M$. If $\widetilde M$ is of finite
type at
$0$ and $H$ is not totally degenerate at $0$, then $H$ is CR
transversal to $\widetilde M$ at $0$. 
\end{Thm}

An inspection of the proof of Theorem \ref{t:main} shows
that Theorem \ref{t:hyper} is a consequence of the
following stronger version of Lemma \ref{l:B} in the
hypersurface case, which was originally proved in
\cite{BRgeom} (Lemma (2.2)). The proof given here is a slight
simplification of that in \cite {BRgeom}.
Since $d=1$, we may assume that the
vector $V\in\bR^d$ in the proof of Theorem \ref{t:main} is
one. 

\begin{Lem}\Label{l:B2} Let $M,\widetilde M$ be formal real
hypersurfaces  through $0\in
\bC^N$, and
$H\colon (\bC^N,0)\to (\bC^N,0)$ a formal holomorphic mapping
with $H(M)\subset\widetilde M$. Assume that $\widetilde M$ is
of finite type at
$0$ and $H$ is not totally degenerate at $0$. If 
\begin{multline}\Label{e:LTid1}\frac{\partial
\ov{\widetilde Q}}{\partial \tilde
\chi}(0,F(z,0),0)\ov{\frac{\partial F}{\partial
w}}(0)\frac{\partial Q}{\partial z}(z,\chi,\bar
Q(\chi,z,0))=\\- \frac{\partial \widetilde Q}{\partial
\tilde z}(F(z,0),\bar F(\chi,\ov Q(\chi,z,0)), \ov
G(\chi,\ov Q(\chi,z,0)))\times \\
\left(\frac{\partial F}{\partial z}(z,0)+\frac{\partial
F}{\partial w}(z,0)\frac{\partial Q}{\partial z}(z,\chi,\ov
Q(\chi,z,0))\right )
\end{multline}
then
\begin{equation}\Label{e:contra31} \frac{\partial \widetilde
Q}{\partial \tilde z}(0,\ov{F}(\chi,0),0)\equiv 0.
\end{equation}
\end{Lem}

\begin{proof} Let us denote by $C(z,\chi)$ the $n\times n$
matrix $\displaystyle\frac{\partial F}{\partial
z}(z,0)+\displaystyle\frac{\partial F}{\partial
w}(z,0)\displaystyle\frac{\partial Q}{\partial z}(z,\chi,\ov
Q(\chi,z,0))$ and by $\delta(z,\chi)$ the determinant $\det
C(z,\chi)$. Recall that $\delta(z,0)\not\equiv 0$ since
$H$ is Segre finite. 
Let $E(z,\chi)$ be the $n\times n$ matrix
satisfying the identity
\begin{equation}\Label{e:EC}
C(z,\chi)E(z,\chi)=E(z,\chi)C(z,\chi)=\delta (z,\chi) I,
\end{equation}
where $I$ is the $n\times n$ identity matrix. Then 
\begin{equation}
\frac{\partial Q}{\partial
z}(z,\chi,\bar Q(\chi,z,0))E(z,\chi)C(z,\chi)=
\delta(z,\chi)\frac{\partial Q}{\partial
z}(z,\chi,\bar Q(\chi,z,0)).
\end{equation}
By Cramer's rule, it follows that 
\begin{equation}\Label{e:QzE}
\frac{\partial Q}{\partial
z}(z,\chi,\bar Q(\chi,z,0))E(z,\chi)=\Delta(z,\chi),
\end{equation}
where $\Delta(z,\chi)=(\Delta_1(z,\chi),\ldots,
\Delta_n(z,\chi))$ with $\Delta_i(z,\chi)=\det B_i(z,\chi)$
and $B_i(z,\chi)$ is the $n\times n$ matrix obtained by
replacing the $i$:th row in $C(z,\chi)$ by
$\displaystyle\frac{\partial Q}{\partial z}(z,\chi,\bar
Q(\chi,z,0))$. In view of \eqref{e:contra0}, we may assume,
after a linear transformation applied to the components of
$F$ if necessary, that $\displaystyle \frac{\partial
F}{\partial w}(0)=(1,0,\ldots,0)^t$. We set
\begin{equation}\Label{e:A}
A(z):= \frac{\partial
\ov{\widetilde Q}}{\partial \tilde
\chi}(0,F(z,0),0)\ov{\frac{\partial F}{\partial
w}}(0)=\frac{\partial
\ov{\widetilde Q}}{\partial \tilde
\chi_1}(0,F(z,0),0)
\end{equation}
If we multiply \eqref{e:LTid1} on the right by $E(z,\chi)$
and use \eqref{e:EC} and \eqref{e:QzE}, we obtain
\begin{equation}\Label{e:newLTid}
A(z)\Delta_i(z,\chi)=
- \frac{\partial \widetilde Q}{\partial
\tilde z_i}(F(z,0),\bar F(\chi,\ov Q(\chi,z,0)), \ov
G(\chi,\ov Q(\chi,z,0)))\delta(z,\chi)
\end{equation}
for $i=1,\ldots, n$.
 We shall also make use of
the identity 
\begin{equation}\Label{e:detid}
\delta(z,\chi)=\det \frac{\partial F}{\partial
z}(z,0)+\sum_{i=1}^n\frac{\partial F}{\partial
w}(z,0)\Delta_i(z,\chi),
\end{equation}
which is a consequence of the following general lemma
in linear algebra. 

\begin{Lem}\Label{l:gla} Let $\mathcal C$ be an $n\times n$
matrix and
$x,y$ column vectors in $\bC^n$. Then, 
$$
\det (\mathcal C+xy^t)=\det \mathcal C +\sum_{i=1}^n x_i \det
\mathcal C_i(y),
$$
where $\mathcal C_i(y)$ is the $n\times n$ matrix obtained
from
$\mathcal C$ by replacing the $i$:th row by $y^t$.
\end{Lem}

The proof of Lemma \ref{l:gla} is obtained by cofactor
expansion along the rows. The details are left to the reader.

To finish the proof of Lemma \ref{l:B2}, we shall show that 
$A(z)\equiv 0$, where $A(z)$ is given by \eqref{e:A}. The
conclusion of the lemma will then follow from the
identity \eqref{e:newLTid} and the fact that
$\delta(z,\chi)\not\equiv 0$. Suppose, in order to reach a
contradiction, that $A(z)\not\equiv 0$. We may then choose a
complex analytic curve $t\mapsto \gamma(t)$ through $0\in
\bC^n$ such that $A(\gamma(t))\not\equiv 0$ and
$\delta(\gamma(t),\chi)\not\equiv 0$. We Taylor expand
$A(\gamma(t))$ and $\delta(\gamma(t),\chi)$ in $t$ to obtain
\begin{equation}\Label{e:Adelta}
A(\gamma(t))=at^p+O(t^{p+1}),\quad
\delta(\gamma(t),\chi)=e(\chi)t^k+O(t^{k+1}),
\end{equation}
where $a\neq 0$ and $e(\chi)\not\equiv 0$.  Also, by the identity 
\begin{equation}\Label{e:nothing}
\frac{\partial
\ov{\widetilde Q}}{\partial \tilde
\chi}(0,F(z,0),0)=\ov{\frac{\partial
{\widetilde Q}}{\partial \tilde
z}(0,\ov{F(z,0)},0)},
\end{equation}
and our assumption $A(z)\not\equiv 0$,
we conclude 
\begin{equation}\Label{e:qz1}
\frac{\partial
{\widetilde Q}}{\partial \tilde
z}(0,\bar F(\chi,0),0)\not\equiv 0.
\end{equation}
Let us expand
$\Delta_1(\gamma(t),\chi)$ in
$t$
\begin{equation}
\Delta_1(\gamma(t),\chi)=c(\chi)t^{l}+O(t^{l+1}).
\end{equation}
It follows from the identity \eqref{e:newLTid}, with
$z=\gamma(t)$ and $i=1$, and \eqref{e:qz1} that $c(\chi)\not\equiv 0$ and $p+l=k$.
We conclude that the right hand side of 
\eqref{e:newLTid}, with
$z=\gamma(t)$ and $i=2,\ldots,n$, is $O(t^{p+l})$ and,
hence,  $\Delta_i(\gamma(t),\chi)=O(t^l)$
for  $i=1,\ldots, n$. Since $\displaystyle\frac{\partial
F}{\partial w}(0)=(1,0,\ldots,0)^t$, we obtain
\begin{equation}\Label{e:Fwdelta}
\sum_{i=1}^n\frac{\partial F_i}{\partial
w}(\gamma(t),0)\Delta_i(\gamma(t),\chi)=c(\chi)t^l+
O(t^{l+1}).
\end{equation}
Let $\delta_0(z):=\displaystyle\frac{\partial F}{\partial
z}(z,0)$. Observe that
$\Delta_i(z,0)\equiv 0$ (since $\displaystyle{\frac{\partial
Q}{\partial z}(z,0,\bar Q(0,z,0)}\equiv 0$) and, hence,
$\delta_0(z)=\delta(z,0)$. In particular, 
$\delta_0(\gamma(t))=O(t^k)$ by
\eqref{e:Adelta}. 
It now follows from \eqref{e:detid} and \eqref{e:Fwdelta}
that $l\geq k$. This contradicts the identity $p+l=k$,
established above, since $p\geq 1$, and completes the proof
of Lemma
\ref{l:B2}.
\end{proof}

\section{Notions of
transversality for formal holomorphic mappings}\Label{s:trans}

 In this
section, we shall give some equivalent formulations of CR
transversality for a formal mapping $H\colon (\bC^N,0)\to
(\bC^N,0)$ and a formal generic submanifold $\widetilde M$ through
$0$ in $\bC^N$. We shall prove that CR
tranversality coincides with
transversality  if there is a
formal generic submanifold $M\subset \bC^N$ (of the same
codimension as $\widetilde M$) such that $H(M)\subset \widetilde M$.
Recall that the formal holomorphic mapping $H$ is transversal
to the formal submanifold
$\widetilde M$ at
$0$ if it satisfies \eqref{e:transdef} when regarded as a
real mapping
$(\bR^{2N},0)\to (\bR^{2N},0)$ or, equivalently, if
\begin{equation}\Label{e:transdefH}
\bC T_0\widetilde M + dH(\bC
T_0\bC^N) =
\bC T_0\bC^N.
\end{equation}
Another
equivalent definition for $H$ to be
transversal to
$\widetilde M$ at $0$ is that for any defining
function $\tilde \rho$ of $\widetilde M$ at $0$,
the real formal mapping
$\tilde
\rho\circ H\colon (\bR^{2N},0)\to (\bR^d,0)$ is a submersion
at
$0$. In particular, if $H$ is transversal to $\widetilde M$ at
$0$, then the inverse image
$H^{-1}(\widetilde M)$ is a formal real submanifold of the same
codimension as $\widetilde M$ (and hence of the same dimension). We also note, by
using the standard facts
\begin{equation}
\bC T_0\bC^N=T^{1,0}_0\bC^N\oplus
T^{0,1}_0\bC^N,\quad dH(T^{1,0}_0\bC^N)\subset
T^{1,0}_0\bC^N,
\end{equation}
that if
$H$ is CR transversal to
$\widetilde M$ at $0$, then it is also transversal to $\widetilde M$ at $0$. The converse is not true in general as is seen by
the following  example.

\begin{Ex}\Label{ex:trans} {\rm Let $\widetilde M\subset \bC^2$ be
the real subspace $\bR^2=\{\im \tilde z_1=\im \tilde z_2=0\}$,
which is a generic submanifold, and consider the mapping
$H(z_1,z_2)=(z_1,iz_1)$. Observe that
$H$ is transversal but not CR transversal to $\widetilde M$ at
$0$.
}
\end{Ex}

Note that in Example \ref{ex:trans} the inverse image
$H^{-1}(\widetilde M)$ is the complex subspace $\{z_1=0\}$. Hence,
there is no formal {\it generic} submanifold $M$ of the same
codimension as $\widetilde M$ such
$H(M)\subset \widetilde M$. The following, which is the main result
of this section, shows that transversality is the same as CR
transversality exactly when
$H$ maps a generic submanifold (of the same codimension) into
$\widetilde M$.

\begin{Thm}\Label{t:trans} Let $H\colon (\bC^N,0)\to
(\bC^N,0)$ be a formal holomorphic mapping and $\widetilde M
\subset\bC^N$ a formal generic submanifold of codimension
$d$.  Then the following are equivalent.
\medskip

{\rm (i)} $H$ is CR transversal to $\widetilde M$ at $0$.
\medskip

{\rm (ii)} $H$ is transversal to $\widetilde M$ at $0$ and the
formal real submanifold $H^{-1}(\widetilde M)\subset \bC^N$ is
generic.
\medskip

{\rm (iii)} If $(\tilde z,\tilde w)$ are normal coordinates
for
$\widetilde M$ at $0$ and $H = (F,G)$, then 
$$dG:T^{1,0}_0\bC^N\to
T^{1,0}_0\bC^d$$ has rank $d$.
\medskip

{\rm (iv)} There exists a formal generic submanifold
$M\subset \bC^N$ of codimension $d$ through $0$ such that
$H(M)\subset \widetilde M$ and
\begin{equation}\Label{e:detG}
\det
\displaystyle\frac {\partial G} {\partial w} (0)
\not= 0,
\end{equation}
where
$(z,w)$ are normal coordinates for $M$ at $0$, $(\tilde
z,\tilde w)$ are normal coordinates for $\widetilde M$ at $0$,
and $H(z,w)=(F(z,w),G(z,w))$.
\end{Thm}

\begin{proof} We begin by showing (i)$\iff$(iii). Let
$Z=(Z_1,\ldots, Z_N)$ be coordinates on the source space
$\bC^N$ and $\tilde Z=(\tilde z,\tilde w)$ normal coordinates for
$\widetilde M$ at $0$. By taking $\partial/\partial
Z_1,\ldots,\partial /\partial Z_N$ as a basis for the source
$T^{1,0}_0\bC^N$ and $\partial/\partial \tilde z_1,\ldots,
\partial/\partial \tilde z_n,\partial/\partial \tilde w_1,
\ldots,
\partial/\partial \tilde w_d$ as a basis for the target
$T^{1,0}_0\bC^N$, we may represent the induced linear mapping
$dH\colon \bC^N\to \bC^n\times \bC^d$ by the matrix
\begin{equation}
dH=\left(\frac{\partial F}
{\partial
Z}(0),
\frac{\partial G}{\partial Z}(0)\right)^t
\end{equation}
The equivalence (i) and (iii) now follows directly from the
fact that $T^{1,0}_0\widetilde M$, in the basis $\partial/\partial \tilde z_1,\ldots,
\partial/\partial \tilde z_n,\partial/\partial \tilde w_1,
\ldots,
\partial/\partial \tilde w_d$, equals the subspace
$\bC^n\times\{0\}\subset \bC^n\times\bC^d\cong
T^{1,0}_0\bC^N$.

Next, we show (i)$\implies$(ii). We already observed in the introduction  that CR
tranversality of $H$ implies transversality and that the
inverse image $M:=H^{-1}(\widetilde M)\subset \bC^N$ is a formal
real submanifold of codimension $d$. Thus, we need only show
that $M$ is generic. For this, we let $(\tilde z,\tilde w)$
be normal coordinates for $\widetilde M$ at $0$ and write
$H=(F,G)$. In these cordinates, $\widetilde M$ is given by the
defining equation
$$
\tilde \rho(\tilde z,\tilde w, {\tilde \chi},{\tilde
\tau}):=\tilde w-\widetilde Q(\tilde z,{\tilde \chi},{\tilde
\tau})=0,
$$ where $\widetilde Q$ satisfies the normalization in
\eqref{e:normal}. The defining function for $M$ is then
$$
\rho(Z,\zeta):=G(Z)-\widetilde Q(F(Z),\bar F(\zeta),\bar
G(\zeta)).
$$
Hence, we have $\partial_Z\rho(0)=\partial_Z G(0)$, and
the fact that $M$ is generic follows from (iii), which has
already been shown to be equivalent to (i).

We now prove (iii)$\iff$(iv). Note that we have shown
(iii)$\implies$(ii), so that in both statements there is a
formal generic submanifold $M\subset \bC^N$ of dimension $d$
that satisfies $H(M)\subset \widetilde M$. Let $Z=(z,w)$ be normal
coordinates for $M$ at $0$, $\tilde Z=(\tilde z,\tilde w)$ normal
coordinates for $\widetilde M$ at $0$, and
$H(z,w)=(F(z,w),G(z,w))$. We take the same bases for the
source and target copies of 
$T^{1,0}_0\bC^N$ as in the proof of (i)$\iff$(iii) except that
 $\partial/\partial Z=(\partial/\partial z,\partial/\partial w$)
 on the source space. The implication
(iv)$\implies$(iii) is immediate. The converse follows from
the fact that $G(z,0)\equiv 0$ and, hence,
$\displaystyle{\frac{\partial G}{\partial z}}(0)=0$.

To complete the proof, we must show (ii)$\implies$ any of the
equivalent (i), (iii), (iv). We use the normal coordinates
$(z,w)$ for
$M:=H^{-1}(\widetilde M)$ and
$(\tilde z,\tilde w)$ for $\widetilde M$ as above and write
$H=(F,G)$. The complexified tangent spaces for $M$ and
$\widetilde M$ can then be written
\begin{equation}
\bC T_0 M=T^{1,0}_0M+T^{0,1}_0M+V_0, \quad \bC T_0
\widetilde M=T^{1,0}_0\widetilde M+T^{0,1}_0\widetilde M+\tilde V_0,
\end{equation}
where $V_0$ and $\tilde V_0$ are the complex spaces
spanned by $\partial/\partial w_1+\partial/\partial \bar
w_1,\ldots, \partial/\partial w_d+\partial/\partial \bar
w_d$ and $\partial/\partial \tilde w_1+\partial/\partial \bar
{\tilde w}_1,\ldots, \partial/\partial
\tilde w_d+\partial/\partial
\bar {\tilde w}_d$, respectively. Let $J$ be
the standard complex structure on $\bC T_0\bC^N$. Then,
since
$M$ and
$\widetilde M$ are generic (see e.g.\ \cite{BER99a}, Proposition
1.3.19 (iii)), we have
\begin{equation} \Label{e:spacedecomp}
\bC T_0\bC^N =\bC T_0M+JV_0,\quad
\bC T_0\bC^N =\bC T_0\widetilde M+J\tilde V_0,
\end{equation}
where both sums are direct. Since $H$ sends $M$ into
$\widetilde M$, the definition of transversality
\eqref{e:transdefH} implies that
\begin{equation}\Label{e:trans1}
\bC T_0\widetilde M+dH(JV_0)=\bC T_0\bC^N.
\end{equation}
Note that $JV_0$ is the complex space spanned by
$\partial/\partial w_1-\partial/\partial \bar
w_1,\ldots, \partial/\partial w_d-\partial/\partial \bar
w_d$ and, hence, $dH(JV_0)$ is spanned by
\begin{multline}
dH\left(\frac{\partial}{\partial w_j}-\frac{\partial}{\partial
\bar w_j}\right)=\\ \sum_{k=1}^n\left( \frac{\partial
F_k}{\partial w_j}(0)\frac{\partial}{\partial \tilde z_k} -
\overline{\frac{\partial F_k}{\partial
w_j}}(0)\frac{\partial}{\partial
\bar {\tilde z}_k} \right)+  \sum_{l=1}^d\left(
\frac{\partial G_l}{\partial w_j}(0)\frac{\partial}{\partial
\tilde w_l} -
\overline{\frac{\partial G_l}{\partial
w_j}}(0)\frac{\partial}{\partial
\bar {\tilde w}_l} \right),
\end{multline}
where $j=1,\ldots d$. As remarked above, $T^{1,0}_0\widetilde M$
and
$T^{0,1}_0\widetilde M$ are spanned by $\partial/\partial \tilde
z_k,\partial/\partial\bar{\tilde
z}_k$, $k=1,\ldots,n$. Thus, using also the fact that
$\displaystyle{\frac {\partial G}{\partial w}(0)}$ is a
real matrix,  we have
\begin{equation}
dH\left(\frac{\partial}{\partial w_j}-\frac{\partial}{\partial
\bar w_j}\right)=  \sum_{l=1}^d
\frac{\partial G_l}{\partial
w_j}(0)\left( \frac{\partial}{\partial
\tilde w_l} -
\frac{\partial}{\partial
\bar {\tilde w}_l} \right)\mod T^{1,0}_0\widetilde
M+T^{0,1}_0\widetilde M,
\end{equation}
for $j=1,\ldots, d$, and hence
$dH(JV_0)\subset J\tilde
V_0$ mod $\bC T_0\widetilde M$ with
equality if and only if
\eqref{e:detG} is satisfied. In view of
the second equality of \eqref{e:spacedecomp} and
\eqref{e:trans1}, we see that (i)$\implies$(iv), which
completes the proof of Theorem
\ref{t:trans}.
\end{proof}

We conclude this section by relating the notion of CR
transversality, given in the introduction \eqref{e:CRtrans0},
for a smooth CR mapping
$f\colon (M,0)\to (\bC^N,0)$ with $f(M)\subset \widetilde M$,
where
$M,\widetilde M\subset \bC^N$ are smooth generic submanifolds of
the same dimension, to that of the induced formal mapping
$H\colon (\bC^N,0)\to (\bC^N,0)$ (see Section \ref{s:prel1}).
We have the following.

\begin{Thm}\Label{t:CRtrans} Let $M,\widetilde M\subset \bC^N$ be
smooth generic submanifolds through $0$ of the same codimension.
Suppose
$f\colon (M,0)\to (\bC^N,0)$ is a smooth CR mapping with
$f(M)\subset
\widetilde M$ and $H\colon (\bC^N,0)\to (\bC^N,0)$ the
induced formal mapping. The following are equivalent.
\medskip

{\rm (i)} $f$ is CR transversal to $\widetilde M$ at $0$.
\smallskip

{\rm (ii)} $H$ is CR transversal to $\widetilde M$ at $0$.
\smallskip

{\rm (iii)} $H$ is transversal to $\widetilde M$ at $0$.
\end{Thm}

The equivalence of (ii) and (iii) follows from Theorem \ref
{t:trans}.  The rest of the proof of Theorem
\ref{t:CRtrans} is left to the reader.

\begin{Rem}\Label{r:diffdim} {\rm The notion of CR
transversality can be generalized to formal mappings between
complex spaces of different dimensions. If $H\colon
(\bC^k,0)\to (\bC^N,0)$ is a formal holomorphic mapping and
$\widetilde M$ is a formal generic submanifold through $0$ in
$\bC^N$, we say that $H$ is CR transversal to $\widetilde M$ at $0$
if $T^{1,0}_0\widetilde M+dH(T^{1,0}_0\bC^k)=T^{1,0}_0\bC^N$. Note,
as before, that if $H$ is transversal to $\widetilde M$ at $0$, then
$H^{-1}(\widetilde M)\subset \bC^k$ is a formal real
submanifold of the same codimension as $\widetilde M$. The
statement and proof of Theorem \ref{t:trans} carry over to
this case also. Similarly,
the notion of CR transversality for a smooth CR mapping
$f\colon (M,0)\to (\bC^N,0)$ sending $M$ into $\widetilde M$
given by
\eqref{e:CRtrans0} can be extended to the case where $M$ is a
smooth generic submanifold of $\bC^k$, where $k$ is not
necessarily equal to $N$. The statement of Theorem
\ref{t:CRtrans} also applies in this case. 
 }
\end{Rem}

\section {Applications to nondegeneracy of mappings between
essentially finite manifolds}\Label{s:new}

In this section, we prove some geometric properties of
sufficiently nondegenerate formal mappings
sending a generic essentially finite manifold into a
generic manifold of the same dimension. We begin by recalling 
the definition of essential finiteness of a formal generic
submanifold $M$ through $0$ in $\bC^N$. We choose normal
coordinates $(z,w)\in \bC^n\times \bC^d$ for $M$ at
$0$ as described in Section \ref{s:prel1}, so that $M$ is
defined by
\eqref{e:def}, and expand the power series
$Q(z,\chi,\tau)=(Q^1(z,\chi,\tau), \ldots,
Q^d(z,\chi,\tau))^t$ at
$\tau=0$ as a Taylor series in $z=(z_1,\ldots,z_n)$:
\begin{equation}\Label{e:expQ} Q^j(z,\chi,0)=\sum_{\alpha\in
\mathbb N_+^n}q^j_\alpha(\chi)z^\alpha,\quad j=1,\ldots, d.
\end{equation} We also write
$q_\alpha(\chi)=(q^1_\alpha(\chi),\ldots,
q^d_\alpha(\chi))^t$. Observe that all the power series
$q^j_\alpha(\chi)$ have zero constant terms. If we
consider $\bC\dbl \chi\dbr$ as a ring, then the {\it essential
ideal} $I_M$ of $M$ at $0$ (see e.g.\ \cite{BER99a}) is the
ideal generated by all the $q_\alpha^j(\chi)$. This ideal is
an invariant of $M$ (cf.\ \cite{BER99a}; see also
Proposition
\ref{p:Essinv} below), and $M$ is said to be {\it essentially
finite} at $0$ if $I_M$ has finite codimension in $\bC\dbl
\chi\dbr$. Its codimension is then called the {\it essential
type} of $M$ at $0$ and is denoted $\ess_0(M)$.

If 
$H\colon (\bC^N,0)\to (\bC^N,0)$ is a formal finite mapping,
then we define the {\it multiplicity}
of $H$ by 
\begin{equation}\Label{e:mH} \mult(H):=\dim_\bC
\bC\dbl Z\dbr/I(H(Z)).
\end{equation}
 For a holomorphic mapping $H$, $\mult(H)$ is the
number of distinct preimages of a generic point.

Before we state our main result, we introduce
another, weaker, transversality notion for a formal mapping. Let
$\widetilde M$ be a formal generic submanifold of codimension $d$
through $0$ in $\bC^N$, and let
$(\tilde z,\tilde w)$ be normal coordinates for
$\widetilde M$ at $0$. If $H\colon (\bC^N,0)\to
(\bC^N,0)$ is a formal mapping, with $H=(F,G)$ in the normal
coordinates $(\tilde z,\tilde w)$, then we shall say
that
$H$ is {\it transversally regular} to $\widetilde M$ at $0$ if the
ideal
$I(G(Z))$ generated by
$G_1(Z),\ldots, G_d(Z)$ has dimension $n=N-d$, i.e.\ if 
 the Krull dimension of the ring
$
\bC\dbl Z\dbr/I(G(Z)) $ is $N-d$. An equivalent way of
expressing this is that the ideal
$I(\tilde
\rho(H(Z),0))$ generated by
$\tilde \rho^1(H(Z),0),\ldots, \tilde \rho^d(H(Z),0)$ has
dimension $N-d$; here,
$\tilde
\rho=(\tilde \rho^1,\ldots, \tilde \rho^d)$ is any defining
function for $\widetilde M$, and $Z$ is any local coordinate in
$\bC^N$ near $0$. Another equivalent condition is that 
$\{G_1(Z),\ldots, G_d(Z)\}$ is a regular sequence 
in $\bC\dbl Z\dbr/I(H(Z))$.

  We should
point out that if
$H$ is CR transversal to
$\widetilde M$ at $0$, then it is also transversally regular to
$\widetilde M$ at $0$. 
For a convergent mapping $H$, being
transversally regular to $\widetilde M$ at $0$ is also equivalent
to the dimension of the germ at $0$ of the variety
$\{G_1(Z)=\ldots=G_d(Z)=0\}$ (or equivalently the variety
$\{\tilde\rho^1(H(Z),0)=\ldots=\tilde
\rho^d(H(Z),0)=0\}$) being $n$. If $\widetilde M$ is a
formal hypersurface (i.e.\ $d=1$), then $H$ is transversally
regular to $\widetilde M$ at $0$ if and only if $G\not\equiv0$.
The main theorem of
this section is the following.

\begin{Thm}\Label{t:Hfinite} Let $M,\widetilde M$ be formal
generic submanifolds of the same CR dimension through $0\in
\bC^N$, and
$H\colon (\bC^N,0)\to (\bC^N,0)$ a formal holomorphic mapping
sending $M$ into $\widetilde M$. Assume that $M$ is essentially
finite and $\widetilde M$ is of finite type at $0$. If $H$
is transversally regular to $\widetilde M$ at $0$, then
$H$ is a finite mapping that is CR transversal to $\widetilde M$ at
$0$. Furthermore, 
$\widetilde M$ is essentially finite at
$0$ with
\begin{equation}\Label{e:MHtM}
\ess_0(M)=\mult(H)\ess_0(\widetilde M).
\end{equation}  In particular, $\ess_0(M)=\ess_0\widetilde M$ if
and only if
$H$ is a formal biholomorphism.
\end{Thm}

For the proof of Theorem \ref{t:Hfinite}, we need some
preliminary results. We begin by taking normal coordinates
$Z=(z,w)$ for $M$ and $\tilde Z=(\tilde z,\tilde w)$ for
$\widetilde M$, and letting $H$ be a formal mapping sending
$M$ into
$\widetilde M$. If $H$ is Segre finite, then we denote by
$m_H$ the integer \begin{equation}\Label{e:mH} m_H:=\dim_\bC
\bC\dbl z\dbr/I(F(z,0)).
\end{equation}
It should be noted that the definition of 
$m_H$ is independent of the choices of normal
coordinates for $M$ and $\widetilde M$. Observe that if
$M$ and
$\widetilde M$ are real-analytic and
$H$ is a holomorphic mapping, then the Segre varieties
$\Sigma_0$, $\tilde \Sigma_0$ of $M$, $\widetilde M$ are given by
$\Sigma_0=\{(z,0)\}$, $\tilde \Sigma_0=\{(\tilde z,0)\}$,
respectively, and the induced mapping
of Segre varieties $h\colon (\Sigma_0,0)\to (\tilde
\Sigma_0,0)$ is given by $h(z)=F(z,0)$. Hence, $H$ is 
Segre finite if and only if $h$ is a finite mapping.

The following useful observation is standard.

\begin{Lem}\Label{l:wellknown} Let $M,\widetilde M$ be formal
generic submanifolds of the same CR dimension through $0\in
\bC^N$, and
$H\colon (\bC^N,0)\to (\bC^N,0)$ a formal holomorphic mapping
sending $M$ into $\widetilde M$. The following hold.
\medskip

{\rm (i)} If $H$ is finite, then $H$ is Segre finite,  
$m_H\leq\mult(H)$, and $H$ is transversally regular to
$\widetilde M$ at $0$.
\smallskip

{\rm (ii)} If $H$ is Segre finite, then $H$ is
not totally
degenerate.

{\rm (iii)} If $H$ is Segre finite and CR
transversal to $\widetilde M$ at $0$, then $H$ is finite and
$m_H=\mult(H)$.
\end{Lem}

\begin{proof} We first prove (i). Since, as is well known,
$G(z,w)=a(z,w)w$ for some $d\times d$ matrix valued formal
power series
$a(z,w)$, we have
$I(G(z,w))\subset I(w)$. It follows that
$I(F(z,w),G(z,w))\subset I(F(z,0),w)$. From this,
the fact that $H$ is Segre finite  and
$m_H\leq \mult(H)$ immediately follows. Since
$I(F(z,w),G(z,w))$  being of finite codimension is
equivalent to having dimension $0$, it follows, in
particular, that $\{G_1(Z),\ldots, G_d(Z)\}$ is a regular
sequence, i.e.\ $H$ is transversally regular.

The statement (ii) is an immediate consequence of the fact
that a finite map $h$ satisfies
$\Jac(h)\not\equiv 0$; see e.g.\ \cite{BER99a}, Theorem
5.1.37.

For (iii), we observe that the CR transversality of
$H$ implies that the ideal $I(G(z,w))$ equals
$I(w)$. Hence, $I(F(z,w),G(z,w))=I(F(z,0),w)$,
from which (iii) easily follows.
\end{proof}

We shall also need the following result relating the
essential types of $M$ and $\widetilde M$ when there is a formal
holomorphic map
$H$ from $M$ to $\widetilde M$. This result is essentially due to
Baouendi and the second author \cite{BRgeom} in
the case where $M$ and
$\widetilde M$ are hypersurfaces. See also Meylan 
\cite{Me1} for manifolds of higher codimension.

\begin{Pro}\Label{p:Essinv} Let $M,\widetilde M$ be formal generic
submanifolds of the same CR dimension through $0\in \bC^N$,
and
$H\colon (\bC^N,0)\to (\bC^N,0)$ a formal holomorphic mapping
sending $M$ into $\widetilde M$. Assume that $H$ is CR
transversal to $\widetilde M$ at
$0$. If $M$ is essentially finite at $0$, then so is $\widetilde M$ and $H$ is Segre finite. If $\widetilde M$ is
essentially finite at $0$ and $H$ is Segre finite,
then $M$ is essentially finite at $0$. In both cases, we have
the identity
\begin{equation}\Label{e:Essid2}
\ess_p(M)=m_H\ess_p(\widetilde M),
\end{equation} where $m_H$ is defined by \eqref{e:mH}.
\end{Pro}

\begin{proof} As before, we take normal coordinates
$Z=(z,w)$ and $\tilde Z = (\tilde z,\tilde w)$ for $M$
and $\widetilde M$ respectively.  Since
$G(z,0)\equiv 0$ and
$\det(\partial G/\partial w)(0)\neq 0$, we have
\begin{equation}\Label{e:Gexp} G(z,w)=a(z,w)w,
\end{equation} where the $d\times d$ matrix
$a(z,w)$ satifies $\det\, a(0)\neq 0$.
Moreover, since $\partial
\widetilde Q/\partial \tilde z(0)=0$, we also have
\begin{equation}\Label{e:QFexp}
\widetilde Q(F(z,w),\tilde\chi,0)=\widetilde Q(F(z,0),\tilde
\chi,0)+b(z,\tilde \chi,w)w,
\end{equation} where the $d\times d$ matrix
$b(z,\tilde\chi,w)$ satisfies
$b(0)=0$. Substituting \eqref{e:Gexp} and
\eqref{e:QFexp} with $w=Q(z,\chi,0)$ and $\tilde \chi=\bar
F(\chi,0)$ in the identity \eqref{e:basic1} above, we obtain
\begin{multline}\Label{e:refsub}
 a(z,Q(z,\chi,0))Q(z,\chi,0)=
 \\\widetilde Q(F(z,0),\bar
 F(\chi,0),0)+b(z,\bar F(\chi,0),Q(z,\chi,0))Q(z,\chi,0),
\end{multline} which can be written as
\begin{equation}\Label{e:refsub1}
 \widetilde Q(F(z,0),\bar
 F(\chi,0),0)=c(z,\chi)Q(z,\chi,0)
\end{equation} where $c(z,\chi)$ is a $d\times d$ matrix with
$\det\, c(0)\neq 0$. Now, the conclusions of 
Proposition \ref{p:Essinv} follow from
\eqref{e:refsub1} by the algebraic arguments in Section 4 of
\cite{BRgerms}. We 
refer the reader to that paper for the details.
\end{proof}

The final result we need for the proof of Theorem
\ref{t:Hfinite} is the following, which can be viewed as a
generalization of a result in
\cite{BRgeom} to higher codimensional manifolds.

\begin{Thm}\Label{t:transfinite} Let $M,\widetilde M$ be formal
generic submanifolds of the same dimension through $0\in
\bC^N$, and
$H\colon (\bC^N,0)\to (\bC^N,0)$ a formal holomorphic mapping
sending $M$ into $\widetilde M$. If $M$ is essentially finite at
$0$ and $H$ is transversally regular, then $H$ is Segre finite.
\end{Thm}

\begin{proof} Assume that $H$ is not Segre finite,
i.e.\ that the ideal $I(F(z,0))$, or equivalently $I(\bar
F(\chi,0))$, does not have finite codimension. We will show
that $H$ is not transversally regular. The ideal
$I(\bar F(\chi,0))$ does not have finite codimension in
$\bC\dbl \chi\dbr$ if and only if there is a (not identically
zero) formal holomorphic mapping $\mu\colon (\bC,0)\to
(\bC^n,0)$ such that
$\bar F(\mu(t),0)\equiv 0$ (see e.g.\ Lemma 3.32 in
\cite{BER00}). Hence, since
$I(\bar F(\chi,0))$ is assumed not to have finite codimension,
we conclude that there is such a mapping $\mu(t)$. If we
substitute $\chi=\mu(t)$ into the equation \eqref{e:basic1},
we conclude (since $\widetilde Q(z,0,0)\equiv 0$) that
\begin{equation}\Label{e:Gis0} G(z,Q(z,\mu(t),0))\equiv 0.
\end{equation} Now, we claim that the Krull dimension of
$\bC\dbl z,w\dbr/I(G(z,w))$ is $>n$. Since $G(z,0)\equiv 0$,
we have that
$I(G(z,w))\subset I(w)$, which is a prime ideal. Since
$\bC\dbl z,w\dbr/I(w)\cong \bC\dbl z\dbr$,  the
Krull dimension of $\bC\dbl z,w\dbr/I(w)$ is equal to $n$.
Thus, to show that the Krull dimension of $\bC\dbl
z,w\dbr/I(G(z,w))$ is
$>n$, we must show that there is a prime ideal $\frak p$ such
that
$I(G(z,w))\subset \frak p\subset I(w)$ and $\frak p\neq
I(w)$. Let
$R\colon (\bC^{n+1},0)\to (\bC^N,0)$ be the formal holomorphic
mapping $R(z,t)=(z,Q(z,\mu(t),0))$ and  $\psi\colon \bC\dbl
z,w\dbr\to
\bC\dbl z,t\dbr$  be
the ring homomorphism defined by $\psi(h) = h\circ
R$. Put
$\frak p:=
\ker \psi$. Then, $\frak p$ is a prime ideal and
$I(G(z,w))\subset \frak p$ by \eqref{e:Gis0}. Moreover, since
$R(z,0)=(z,Q(z,0,0))=(z,0)$, it follows that $\frak p\subset
I(w)$. To complete the proof, we must show that $\frak p\neq
I(w)$. To see this, observe that by
definition $w_i\in \frak p$ if and only if
$Q^i(z,\mu(t),0)\equiv 0$. But, the latter is equivalent to
saying that
$q^i_\alpha(\mu(t))=0$, for all multi-indices $\alpha$. This
cannot happen, for every $i=1,\ldots d$, since $M$ is
essentially finite (which, as the reader should recall, means
that the ideal generated by all the
$q^i_\alpha(\chi)$ has finite codimension in $\bC\dbl
\chi\dbr$). This completes the proof.
\end{proof}

\begin{proof}[Proof of Theorem $\ref{t:Hfinite}$] If $H$ is
transversally regular to $\widetilde M$ at
$0$, then it is also Segre finite by Theorem
\ref{t:transfinite}. It follows from Theorem
\ref{t:main} that $H$ is CR transversal, and the conclusion of
Theorem \ref{t:Hfinite} follows from Proposition
\ref{p:Essinv} and Lemma \ref{l:wellknown} (iii).\end{proof}

We would like to point out that there is no assumption of
finite type on either $M$ or $\widetilde M$ in Proposition
\ref{p:Essinv} or in Theorem
\ref{t:transfinite}.

By using Theorem \ref{t:Hfinite}, we can give a sufficient
condition on
$M$ that guarantees that every finite formal mapping sending $M$
into another generic submanifold of the same dimension is a
formal biholomorphism.

\begin{Thm}\Label{t:fndsimple} Let $M$ be a formal generic
submanifold through $0\in \bC^N$ and assume that $M$ is
finitely nondegenerate and of finite type at $0$. Then any
formal finite holomorphic mapping $H\colon (\bC^N,0)\to
(\bC^N,0)$ sending $M$ into another formal generic
submanifold $\widetilde M$ (through
$0$) of the same dimension is a local
formal biholomorphism.
\end{Thm}

\begin{proof} Since a finite mapping $H$ satisfies $\Jac
(H)\not\equiv 0$ (see e.g.\ \cite{BER99a}, Theorem 5.1.37), it
follows from Proposition \ref{p:ft}(a) that $\widetilde M$ is of
finite type at $0$. Moreover, Lemma \ref{l:wellknown} (i)
implies that
$H$ is transversally regular . Hence, by Theorem
\ref{t:Hfinite} the multiplicity of $H$ divides the
essential type of $M$ at $0$. Moreover, by Proposition
11.8.27 in \cite{BER99a}, $M$ is essentially finite and
$\ess_0(M)=1$. This proves that
$\mult(H)=1$, which means that $H$ is a formal biholomorphism.
\end{proof}

We note that  the
condition of finite nondegeneracy in Theorem \ref{t:fndsimple}
cannot be replaced by the weaker condition of essential
finiteness. Indeed, consider $H\colon (\bC^2,0)\to (\bC^2,0)$ given by
$H(z,w)=(z^2,w)$, and let $M,\widetilde M\subset
\bC^2$ be given by
\begin{equation}
M\colon \im w=|z|^4,\quad \widetilde M\colon \im w=|z|^2,
\end{equation} Then $H(M) \subset \widetilde M$, but is not a local biholomorphism.  A
more precise condition guaranteeing that the conclusion of Theorem
\ref{t:fndsimple} holds will be given in the forthcoming paper \cite{ER2}.

\section{Closing remarks and proofs of Theorems
\ref{t:maincor}, \ref{t:maincor1}, \ref{t:maincor2},
\ref{t:maincor3}.}\Label{s:last}

For the proofs of Theorems \ref{t:maincor}, \ref{t:maincor1},
\ref{t:maincor2},
\ref{t:maincor3}, we shall assume, without loss of generality, that $p_0=\tilde
p_0=0$. We shall regard $M,
\widetilde M\subset \bC^N$ as formal generic submanifolds and
$H\colon (\bC^N,0)\to(\bC^N,0)$ as a formal holomorphic
mapping, as explained in Section \ref{s:prel1}.

To prove Theorem \ref{t:maincor}, we observe from Lemma
\ref{l:wellknown} (i) that any finite formal mapping
is Segre finite. Theorem \ref{t:maincor} is now a
direct consequence of Theorem \ref{t:main}.
Theorem \ref{t:maincor1} is an immediate consequence of
Theorem
\ref{t:fndsimple}.

For Theorem \ref{t:maincor2}  the
implication (ii)$\implies$(i) 
follows immediately from Theorem
\ref{t:maincor}, as already observed in the introduction. 
 The remaining implication (i)$\implies$(ii)
follows from Proposition \ref{p:Essinv} and Lemma
\ref{l:wellknown} (iii). We should point out that the
assumption that $M$ is of finite type at $0$ is not needed
for (i)$\implies$(ii).
To prove Theorem \ref{t:maincor3}, we observe that 
the statement ``$h$ is finite" means exactly that
$H$ is Segre finite, as noted in Section
\ref{s:new}. Theorem  \ref{t:maincor3} now follows from
Theorem \ref{t:main}.

We conclude this paper with a few remarks and questions. For
the case of a hypersurface ($d=1$), $H=(F,G)$ is
transversally regular (as defined in Section
\ref{s:new}) if and only if   
$G \not\equiv 0$. Hence, for that case, if
$\Jac(H)\not\equiv 0$, then $H$ must be transversally regular to
$\widetilde M$ at $0$. (However, if the codimension of
$\widetilde M$ is $>1$,  a mapping $H\colon(\bC^N,0)\to
(\bC^N,0)$ need not be transversally regular to $\widetilde
M$ even if
$\Jac(H)\not\equiv 0$.)  By Proposition
\ref{p:ft} (a) and  Theorem \ref{t:Hfinite}, we conclude the
following (since an essentially finite hypersurface is necessarily
 also of finite type).
\medskip

(*) {\it If $M\subset\bC^N$ is an
essentially finite hypersurface through $0$ and
$H\colon(\bC^N,0)\to (\bC^N,0)$ satisfies
$\Jac(H)\not\equiv 0$ and sends
$M$ into another hypersurface $\widetilde M$, then $H$ is finite
and CR transversal to $\widetilde M$ at $0$.}\medskip

The above result already follows
from \cite{BRgeom}. One is led to formulate the
following two open questions. The first is whether the
condition that $M$ is essentially finite in (*) can be
weakened to the condition that $M$ is of finite type, i.e.\
\medskip

{\bf Question 1:} {\it Suppose $M\subset\bC^N$ is a hypersurface
of finite type through $0$ and
$H\colon(\bC^N,0)\to (\bC^N,0)$ satisfies
$\Jac(H)\not\equiv 0$ and sends
$M$ into another hypersurface $\widetilde M$. Does it follow that
$H$ is CR transversal?}
\medskip

We note that $H$, in Question 1, need not
be a finite map. For example, let $H:\C^3 \to \bC^3$ be given
(implicitly) by
$H(z_1,z_2,w) = (z_1,z_2w,w)$, and $M,\widetilde M\subset
\bC^3$ the hypersurfaces  given by
\begin{equation}
M\colon \im w=|z_1|^2 +|z_2w|^2,\quad \widetilde M\colon \im
\tilde w=|\tilde z_1|^2+|\tilde z_2|^2,
\end{equation} Then $H(M) \subset \widetilde M$, but $H$ is not finite.  

 The other question concerns the
possible generalization of statement (*).

\medskip

{\bf Question 2:} {\it Suppose $M\subset\bC^N$ is an
essentially finite generic submanifold of finite type through
$0$ and
$H\colon(\bC^N,0)\to (\bC^N,0)$ satisfies
$\Jac(H)\not\equiv 0$ and sends
$M$ into another generic submanifold $\widetilde M\subset \bC^N$
of the same dimension. Does it follow that
$H$ is CR transversal? In particular, if $M$ is assumed to be
finitely nondegenerate, is $H$ necessarily a formal
biholomorphism?}

\end{document}